\documentclass[a4paper, british]{amsart}

\usepackage{libertine}
\usepackage[libertine]{newtxmath}

\usepackage{babel}
\usepackage{enumitem}
\usepackage{hyperref}
\usepackage[utf8]{inputenc}
\usepackage{newunicodechar}
\usepackage{mathtools}
\usepackage{varioref}
\usepackage[arrow,curve,matrix]{xy}

\usepackage{colortbl}
\usepackage{graphicx}
\usepackage{tikz}

\usepackage{mathrsfs}

\definecolor{linkred}{rgb}{0.7,0.2,0.2}
\definecolor{linkblue}{rgb}{0,0.2,0.6}

\numberwithin{figure}{section}

\usepackage[hyperpageref]{backref}

\setdescription{labelindent=\parindent, leftmargin=2\parindent}
\setitemize[1]{labelindent=\parindent, leftmargin=2\parindent}
\setenumerate[1]{labelindent=0cm, leftmargin=*, widest=iiii}

\newunicodechar{α}{\ensuremath{\alpha}}
\newunicodechar{β}{\ensuremath{\beta}}
\newunicodechar{χ}{\ensuremath{\chi}}
\newunicodechar{δ}{\ensuremath{\delta}}
\newunicodechar{ε}{\ensuremath{\varepsilon}}
\newunicodechar{Δ}{\ensuremath{\Delta}}
\newunicodechar{η}{\ensuremath{\eta}}
\newunicodechar{γ}{\ensuremath{\gamma}}
\newunicodechar{Γ}{\ensuremath{\Gamma}}
\newunicodechar{ι}{\ensuremath{\iota}}
\newunicodechar{κ}{\ensuremath{\kappa}}
\newunicodechar{λ}{\ensuremath{\lambda}}
\newunicodechar{Λ}{\ensuremath{\Lambda}}
\newunicodechar{ν}{\ensuremath{\nu}}
\newunicodechar{μ}{\ensuremath{\mu}}
\newunicodechar{ω}{\ensuremath{\omega}}
\newunicodechar{Ω}{\ensuremath{\Omega}}
\newunicodechar{π}{\ensuremath{\pi}}
\newunicodechar{Π}{\ensuremath{\Pi}}
\newunicodechar{φ}{\ensuremath{\phi}}
\newunicodechar{Φ}{\ensuremath{\Phi}}
\newunicodechar{ψ}{\ensuremath{\psi}}
\newunicodechar{Ψ}{\ensuremath{\Psi}}
\newunicodechar{ρ}{\ensuremath{\rho}}
\newunicodechar{σ}{\ensuremath{\sigma}}
\newunicodechar{Σ}{\ensuremath{\Sigma}}
\newunicodechar{τ}{\ensuremath{\tau}}
\newunicodechar{θ}{\ensuremath{\theta}}
\newunicodechar{Θ}{\ensuremath{\Theta}}
\newunicodechar{ξ}{\ensuremath{\xi}}
\newunicodechar{Ξ}{\ensuremath{\Xi}}
\newunicodechar{ζ}{\ensuremath{\zeta}}

\newunicodechar{ℓ}{\ensuremath{\ell}}
\newunicodechar{ï}{\"{\i}}

\newunicodechar{⁰}{\ensuremath{^0}}
\newunicodechar{¹}{\ensuremath{^1}}
\newunicodechar{²}{\ensuremath{^2}}
\newunicodechar{³}{\ensuremath{^3}}
\newunicodechar{⁴}{\ensuremath{^4}}
\newunicodechar{⁵}{\ensuremath{^5}}
\newunicodechar{⁶}{\ensuremath{^6}}
\newunicodechar{⁷}{\ensuremath{^7}}
\newunicodechar{⁸}{\ensuremath{^8}}
\newunicodechar{⁹}{\ensuremath{^9}}
\newunicodechar{ⁱ}{\ensuremath{^i}}

\newunicodechar{⌈}{\ensuremath{\lceil}}
\newunicodechar{⌉}{\ensuremath{\rceil}}
\newunicodechar{⌊}{\ensuremath{\lfloor}}
\newunicodechar{⌋}{\ensuremath{\rfloor}}

\newunicodechar{≅}{\ensuremath{\cong}}
\newunicodechar{⇔}{\ensuremath{\Leftrightarrow}}
\newunicodechar{∃}{\ensuremath{\exists}}
\newunicodechar{±}{\ensuremath{\pm}}

\DeclareFontFamily{OMS}{rsfs}{\skewchar\font'60}
\DeclareFontShape{OMS}{rsfs}{m}{n}{<-5>rsfs5 <5-7>rsfs7 <7->rsfs10 }{}
\DeclareSymbolFont{rsfs}{OMS}{rsfs}{m}{n}
\DeclareSymbolFontAlphabet{\scr}{rsfs}
\DeclareSymbolFontAlphabet{\scr}{rsfs}

\DeclareFontFamily{U}{mathx}{\hyphenchar\font45}
\DeclareFontShape{U}{mathx}{m}{n}{
      <5> <6> <7> <8> <9> <10>
      <10.95> <12> <14.4> <17.28> <20.74> <24.88>
      mathx10
      }{}
\DeclareSymbolFont{mathx}{U}{mathx}{m}{n}
\DeclareFontSubstitution{U}{mathx}{m}{n}
\DeclareMathAccent{\wcheck}{0}{mathx}{"71}

\newunicodechar{∂}{\ensuremath{\partial}}
\newunicodechar{∇}{\ensuremath{\nabla}}

\newunicodechar{↺}{\ensuremath{\circlearrowleft}}
\newunicodechar{∞}{\ensuremath{\infty}}
\newunicodechar{⊕}{\ensuremath{\oplus}}
\newunicodechar{⊗}{\ensuremath{\otimes}}
\newunicodechar{•}{\ensuremath{\bullet}}
\newunicodechar{Λ}{\ensuremath{\wedge}}
\newunicodechar{↪}{\ensuremath{\into}}
\newunicodechar{→}{\ensuremath{\to}}
\newunicodechar{↦}{\ensuremath{\mapsto}}
\newunicodechar{⨯}{\ensuremath{\times}}
\newunicodechar{∪}{\ensuremath{\cup}}
\newunicodechar{∩}{\ensuremath{\cap}}
\newunicodechar{⊋}{\ensuremath{\supsetneq}}
\newunicodechar{⊇}{\ensuremath{\supseteq}}
\newunicodechar{⊃}{\ensuremath{\supset}}
\newunicodechar{⊊}{\ensuremath{\subsetneq}}
\newunicodechar{⊆}{\ensuremath{\subseteq}}
\newunicodechar{⊂}{\ensuremath{\subset}}
\newunicodechar{⊄}{\ensuremath{\not \subset}}
\newunicodechar{≥}{\ensuremath{\geq}}
\newunicodechar{≠}{\ensuremath{\neq}}
\newunicodechar{≫}{\ensuremath{\gg}}
\newunicodechar{≪}{\ensuremath{\ll}}

\newunicodechar{≤}{\ensuremath{\leq}}
\newunicodechar{∈}{\ensuremath{\in}}
\newunicodechar{∉}{\ensuremath{\not \in}}
\newunicodechar{∖}{\ensuremath{\setminus}}
\newunicodechar{◦}{\ensuremath{\circ}}
\newunicodechar{°}{\ensuremath{^\circ}}
\newunicodechar{…}{\ifmmode\mathellipsis\else\textellipsis\fi}
\newunicodechar{·}{\ensuremath{\cdot}}
\newunicodechar{⋯}{\ensuremath{\cdots}}
\newunicodechar{∅}{\ensuremath{\emptyset}}
\newunicodechar{⇒}{\ensuremath{\Rightarrow}}

\newunicodechar{⁰}{\ensuremath{^0}}
\newunicodechar{¹}{\ensuremath{^1}}
\newunicodechar{²}{\ensuremath{^2}}
\newunicodechar{³}{\ensuremath{^3}}
\newunicodechar{⁴}{\ensuremath{^4}}
\newunicodechar{⁵}{\ensuremath{^5}}
\newunicodechar{⁶}{\ensuremath{^6}}
\newunicodechar{⁷}{\ensuremath{^7}}
\newunicodechar{⁸}{\ensuremath{^8}}
\newunicodechar{⁹}{\ensuremath{^9}}
\newunicodechar{ⁱ}{\ensuremath{^i}}

\newunicodechar{⌈}{\ensuremath{\lceil}}
\newunicodechar{⌉}{\ensuremath{\rceil}}
\newunicodechar{⌊}{\ensuremath{\lfloor}}
\newunicodechar{⌋}{\ensuremath{\rfloor}}

\newunicodechar{≅}{\ensuremath{\cong}}
\newunicodechar{⇔}{\ensuremath{\Leftrightarrow}}
\newunicodechar{∃}{\ensuremath{\exists}}
\newunicodechar{±}{\ensuremath{\pm}}

\DeclareFontFamily{OMS}{rsfs}{\skewchar\font'60}
\DeclareFontShape{OMS}{rsfs}{m}{n}{<-5>rsfs5 <5-7>rsfs7 <7->rsfs10 }{}
\DeclareSymbolFont{rsfs}{OMS}{rsfs}{m}{n}
\DeclareSymbolFontAlphabet{\scr}{rsfs}
\DeclareSymbolFontAlphabet{\scr}{rsfs}

\DeclareFontFamily{U}{mathx}{\hyphenchar\font45}
\DeclareFontShape{U}{mathx}{m}{n}{
      <5> <6> <7> <8> <9> <10>
      <10.95> <12> <14.4> <17.28> <20.74> <24.88>
      mathx10
      }{}
\DeclareSymbolFont{mathx}{U}{mathx}{m}{n}
\DeclareFontSubstitution{U}{mathx}{m}{n}
\DeclareMathAccent{\wcheck}{0}{mathx}{"71}

\DeclareMathOperator{\reg}{reg}

\theoremstyle{plain}
\newtheorem{Th}{Théorème}[section]

\newtheorem{cor}[Th]{Corollaire}
\newtheorem{defn}[Th]{Définition}

\newtheorem{lem}[Th]{Lemme}

\newtheorem{Prop}[Th]{Proposition}

\theoremstyle{remark}

\newtheorem{c-n-d}[Th]{Claim and Definition}

\newtheorem{rem}[Th]{\color{blue}{Remarques}}

\newtheorem*{rem-nonumber}{Remark}

\numberwithin{equation}{Th}
\setlist[enumerate]{label=(\thethm.\arabic*), before={\setcounter{enumi}{\value{equation}}}, after={\setcounter{equation}{\value{enumi}}}}
\newcommand{\into}{\hookrightarrow}

\newcommand{\factor}[2]{\left. \raise 2pt\hbox{$#1$} \right/\hskip -2pt\raise -2pt\hbox{$#2$}}

\newcommand{\Publication}[1]{}
\newcommand{\subversionInfo}{}
\newcommand{\svnid}[1]{}
\newcommand{\approvals}[2][Approval]{}
\renewcommand{\phi}{\varphi}
\usepackage{tikz}
\usepackage{tikz-cd}
\tikzset{commutative diagrams/arrow style=Latin Modern}
\author{Mohamed Kaddar} %
\email{\href{mailto:mohamed.kaddar@univ-lorraine.fr}{mohamed.kaddar@univ-lorraine.fr}} %
\title[Sur certains faisceaux de formes méromorphes en géométrie analytique complexe II.]{Sur certains faisceaux de formes méromorphes en géométrie analytique complexe II.}%

\date{\today}
\makeatletter
\hypersetup{
  pdfauthor={\authors},
  pdftitle={\@title},
  pdfsubject={\@subjclass},
  pdfkeywords={\@keywords},
  pdfstartview={Fit},
  pdfpagelayout={TwoColumnRight},
  pdfpagemode={UseOutlines},
  bookmarks,
  colorlinks,
  linkcolor=linkblue,
  citecolor=linkred,
  urlcolor=linkred}
\makeatother
\newcommand{\chapref}[1]{\hyperref[#1]{Chapter~\ref*{#1}}}
\newcommand{\lemmaref}[1]{\hyperref[#1]{Lemme~\ref*{#1}}}
\newcommand{\parref}[1]{\hyperref[#1]{Section~\ref*{#1}}}
\newcommand{\theoremref}[1]{\hyperref[#1]{Théorème~\ref*{#1}}}
\newcommand{\definitionref}[1]{\hyperref[#1]{Définition~\ref*{#1}}}
\newcommand{\propositionref}[1]{\hyperref[#1]{Proposition~\ref*{#1}}}
\newcommand{\conjectureref}[1]{\hyperref[#1]{Conjecture~\ref*{#1}}}
\newcommand{\corollaryref}[1]{\hyperref[#1]{Corollaire~\ref*{#1}}}
\newcommand{\exampleref}[1]{\hyperref[#1]{Exemple~\ref*{#1}}}
\newcommand{\exerciseref}[1]{\hyperref[#1]{Exercise~\ref*{#1}}}
\newcommand{\factref}[1]{\hyperref[#1]{Fact~\ref*{#1}}}
\newcommand{\claimref}[1]{\hyperref[#1]{Claim~\ref*{#1}}}
\newcommand{\remarkref}[1]{\hyperref[#1]{Remark~\ref*{#1}}}
\newcommand{\settingref}[1]{\hyperref[#1]{Setting~\ref*{#1}}}
\newcommand{\appendixref}[1]{\hyperref[#1]{Appendix~\ref*{#1}}}

\theoremstyle{plain}

\theoremstyle{remark}

\newcommand{\eq}[1][r]
   {\ar@<-3pt>@{-}[#1]
    \ar@<-1pt>@{}[#1]|<{}="gauche"
    \ar@<+0pt>@{}[#1]|-{}="milieu"
    \ar@<+1pt>@{}[#1]|>{}="droite"
    \ar@/^2pt/@{-}"gauche";"milieu"
    \ar@/_2pt/@{-}"milieu";"droite"}

\begin{document}

\maketitle
\approvals[Approval for Abstract]{Mohamed & yes}
\begin{abstract}
We study some functorial properties of certain sheaves of meromorphic forms on reduced complex space; particulary, the meromorphic forms which extend holomorphicaly on any désingularisation. The purpose concern their behavior under pull back and higher direct image (and in some case by integration  on fibres of equidimensional or open morphism).
\end{abstract}
\maketitle
\tableofcontents

\noindent
\phantomsection\addcontentsline{toc}{part}{Introduction}
%
%
\vspace{1mm}

\noindent
Dans \cite{K3}, nous avons étudiés certaines propriétés fonctorielles des faisceaux $\omega^{\bullet}_{X}$ sur un espace analytique $X$ réduit de dimension pure $m$ et, principalement, leur comportement par image directe supérieure et par image réciproque. Malheureusement, leur principale lacune, vu la nature dualisante du faisceau $\omega^{m}_{X}$ (rappelons que c'est la $m$-ème homologie  du complexe dualisant et que pour $X$ de Cohen Macaulay, $\omega^{m}_{X}[m]$ est un complexe dualisant), est de ne pas être  stables par image réciproque analytique. D'ailleurs des exemples simples comme la normalisation faible du cusp dans ${\Bbb C}^{2}$ ou la normalisation du "parapluie de Withney" montrent que le relèvement des générateurs de leur faisceau dualisant sont des formes méromorphes à pôles logarithmiques. Néanmoins, on a montré que beaucoup de ces sous faisceaux sont stables par image réciproque de morphismes surjectifs à fibres de dimension constante (équidimensionnels ou ouverts).\vspace{1mm}

\noindent
On se propose, ici, de mener une étude similaire pour le faisceau ${\mathcal L}^{\bullet}_{X}$ sous faisceau de  $\omega^{\bullet}_{X}$ dont les sections se relèvent dans  toute désingularisation en formes méromorphes se prolongeant globalement holomorphiquement. Cette définition de nature globale le rend peu maniable car ne permet pas une description locale à l'instar de $\omega^{\bullet}_{X}$. Il occupe, néanmoins une place importante en géométrie algébrique ou analytique comme on peut le voir dans  les théorèmes d'annulation de Grauert-Riemmenschneider \cite{GR}, dans la théorie des singularités rationnelles, canoniques ou terminales. Les problématiques tournant autour de ces faisceaux  apparaissaient, déjà sous certaines formes, dans l'étude des formes dites de première espèce introduite par Picard  \cite{P.S} et généralisées par Griffiths \cite{Gr}.\vspace{1mm}

 \noindent L'objet de ce qui suit est d'établir un résultat analogue aux {\bf{théorème 1}} et {\bf{théorème 2}} de \cite{K3}, montrant la stabilité par image réciproque arbitraire 
du faisceau  ${\mathcal L}^{\bullet}_{X}$ et son comportement par image directe supérieure. 

Plusieurs approches peuvent être envisagées pour montrer l'existence d'une image réciproque pour ces formes. En effet, si le morphisme $\pi:X\rightarrow S$ est surjectif à fibres de dimension constante, on peut utiliser le {\bf{théorème 1}} de \cite{K3} pour en déduire un pull-back naturel $\pi^{*}{\mathcal L}^{\bullet}_{S}\rightarrow {\mathcal L}^{\bullet}_{X}$. Comme ces faisceaux sont invariants par modifications propres, on en déduit, sans difficultés, grâce aux  théorèmes d'applatissement algébriques ou géométriques, l'existence de cette image réciproque pour tout morphisme propre, surjectif et génériquement équidimensionnel. En utilisant les théorèmes d'applatissement locaux, on peut se passer de la propreté. \vspace{1mm}

\noindent 
Un autre point de vue consiste à proposer une définition et un cadre pour ces formes de sorte à pouvoir récupérer la plupart des propriétés fonctorielles des formes holomorphes habituelles. Pour cela, on est contraint d'imposer des conditions d'incidence aux espaces considérés et de croissance ${\rm L}^{2}$ le long du lieu singulier nous permettant de raisonner génériquement sur des formes holomorphes usuelles; on évite de se retrouver complètement dans le lieu polaire.  On peut, alors, construire cette image réciproque en termes de modifications propres et de désingularisations. \vspace{1mm}

\noindent
Toujours pour un morphisme surjectif mais sans aucune autre condition, on peut utiliser la caractérisation donnée dans \cite{KeSc} des formes holomorphes sur la partie régulière d'un espace complexe réduit dont les images réciproques sur toute désingularisée se prolonge analytiquement et disant que si  $\xi$ est une section de $\Gamma({\rm Reg}(X), \Omega^{k}_{X})$: 
$$\xi\in {\mathcal L}^{k}_{X}\,\Longleftrightarrow\,\xi\wedge \alpha\,\,{\rm et}\, d\xi\wedge \beta\in {\mathcal L}^{m}_{X},\,\forall\,\alpha\in\Omega^{m-k}_{X},\,\,\forall\,\beta \in\Omega^{m-k-1}_{X}$$
Ce théorème est très loin d'être trivial et nécessite des outils très puissants utilisant la théorie des modules de Hodge mixtes.\vspace{1mm}

\noindent
Malheureusement, même si ces  approches donnent déjà un large éventail de morphismes pour lesquelles on a la stabilité désirée par image réciproque, il n'en reste pas moins que le cas d'un morphisme de restriction à un sous espace contenu dans le lieu singulier de l'espace ambiant reste inabordable par ce biais. \vspace{1mm}

\noindent
Dans le cas général d'un  morphisme arbitraire d'espaces complexes réduits, notre construction repose fondamentalement sur la manière dont on va définir une restriction de ces formes méromorphes à pôles dans le lieu singulier de l'espace ambiant  à un sous espace entièrement contenu dans le lieu singulier . L'approche  adoptée est celle qui a été  développée dans \cite{K0}. Elle utilise l'intégration sur les fibres le long du diviseur exceptionnel ou de  la préimage totale. La partie fastidieuse est la vérification de l'indépendance de la construction vis-à-vis de tous les choix que l'on va faire comme la forme de Kähler relative , la désingularisation, les factorisations possibles etc.\vspace{1mm}

\noindent
Tout comme dans \cite{K3}, on établit des résultats d'images directe supérieures pour un morphisme arbitraire d'espaces complexes (réduits) déduits soit de  l'intégration sur les fibres munis de leur structure de cycles (c'est-à-dire dans le cas d'un morphisme géométriquement plat) soit d'une partie émergente de la théorie de la dualité relative. 
On a, alors
\Th{}{}\label{T1}  Soit $\pi:X\rightarrow S$ un morphisme  d'espaces analytiques complexes réduits de dimension pure. Alors,\vspace{1mm}

\noindent
{\bf(i)} Il existe un  unique morphisme de faisceaux 
${\bf{\pi}}^{*}:{\mathcal L}^{\bullet}_{S}\rightarrow \pi_{*}{\mathcal L}^{\bullet}_{X} $ (et même d'algèbre différentielle sur ${\Omega}^{\bullet}_{S}$) prolongeant l'image réciproque des formes holomorphes usuelles et rendant commutatif le diagramme
$$\xymatrix{{\Omega}^{\bullet}_{S}\ar[r]\ar[d]&\pi_{*}{\Omega}^{\bullet}_{X}\ar[d]\\
{\mathcal L}^{\bullet}_{S}\ar[r]&{\pi_{*}}{\mathcal L}^{\bullet}_{X} }$$
De plus, si $\pi$ est surjectif à fibres de dimension  constante, on a aussi, grâce à (\cite{K3}, {\bf{théorème 1}}), le diagramme commutatif
$$\xymatrix{{\mathcal L}^{\bullet}_{S}\ar[r]\ar@{^{(}->}[d] \ar[d]&\pi_{*}{\mathcal L}^{\bullet}_{X}\ar@{^{(}->}[d] \\
{\widetilde{\mathcal L}}^{\bullet}_{S}\ar[r]&\pi_{*}{\widetilde{\mathcal L}}^{\bullet}_{X}}$$
avec ${\widetilde{\mathcal L}}^{\bullet}_{Z}:={\mathcal H}om({\mathcal L}^{{\rm dim}(Z)-\bullet}_{Z}, \omega^{m}_{Z})$ pour $Z=S, X$.
\vspace{1mm}

\noindent {\bf{(ii) Image réciproque et composition.}} 
si $f:X\rightarrow Y$ et $g:Y\rightarrow S$ sont deux morphismes surjectifs d'espaces complexes réduits  alors ${\bf (g\circ f)}^{*}={\bf f}^{*}\circ {\bf g}^{*}$.\vspace{1mm}

\noindent 
\Th{}{}\label{T2} Soit $\pi:X\rightarrow S$ un morphisme  d'espaces analytiques complexes réduits de dimension pure. Alors, si:\vspace{1mm}

\noindent 
{\bf(i)} $\pi$ est  propre génériquement $n$-équidimensionnel, il lui  est associé  un morphisme de faisceaux cohérents $${\mathcal T}^{k}_{\pi, {\mathcal L}}:{\rm I}\!{\rm R}^{n}\pi_{*}{\mathcal L}^{n+k}_{X}\rightarrow{\mathcal L}^{k}_{S} $$
 rendant commutatif le diagramme
$$\xymatrix{ {\rm I}\!{\rm R}^{n}\pi_{*}{\mathcal L}^{n+k}_{X}\ar[r]\ar[d]_{{\mathcal T}^{k}_{\pi, {\mathcal L}}}&{\rm I}\!{\rm R}^{n}\pi_{*}{\omega}^{n+k}_{X}\ar[d]^{{\mathcal T}^{k}_{\pi, {\omega}}}\\
{\mathcal L}^{k}_{S}\ar@{^{(}->}[r]&{\omega}^{k}_{S}}$$ 
(en remplaçant $\pi_{*}$ par $\pi_{!}$ dans le cas non propre) et dont la formation:\vspace{1mm}

\indent
$\bullet$ commute aux restrictions ouvertes sur $X$ (et $S$) et aux changements de base plats,\vspace{1mm}

\indent $\bullet$ est compatible  aux changements de base arbitraires et  à la composition des morphismes d'espaces complexes  dans le sens suivant:\vspace{ 1mm}

\noindent
Pour tout  diagramme commutatif d'espaces analytiques complexes
$$\xymatrix{X_{2}\ar[rr]^{\Psi}\ar[rd]_{\pi_{2}}&&X_{1}\ar[ld]^{\pi_{1}}\\
&S&}$$
 avec $\pi_{1}$ (resp. $\pi_{2}$)  propre, de dimension générique relative $n_{1}$ (resp. $n_{2}$)  et $\Psi$ propre de dimension relative bornée par l'entier   $d:=n_{2}-n_{1}$, il existe un unique morphisme $\tilde\Psi$ rendant commutatif le diagramme  de faisceaux analytiques (cohérents dans le cas propre)
$$\xymatrix{{\rm I}\!{\rm
 R}^{n_{2}}{\pi_{2}}_{!}{\mathcal L}^{n_{2}+q}_{X_{2}}\ar[rr]^{\tilde{\Psi}}\ar[rd]_{{\mathcal T}^{q}_{\pi_{2}}}&&
 {\rm I}\!{\rm
 R}^{n_{1}}{\pi_{1}}_{!}{\mathcal L}^{n_{1}+q}_{X_{1}}\ar[ld]^{{\mathcal T}^{q}_{\pi_{1}}}\\
&{\mathcal L}^{q}_{S}&}$$
En remplaçant $\pi_{!}$ par $\pi_{*}$ dans le cas non propre.
\vspace{1mm}

\noindent
{\bf(ii)} Si $\pi$ est géométriquement plat à fibres de dimension $n$, il est induit par  et induit le morphisme d'intégration sur les fibres.\rm 
\vspace{2mm}

\noindent
\cor{}{}\label{C0} On conserve  les notations et  hypothèses du \theoremref{T2}.
Soient $\Phi$ une famille paracompactifiante de supports de $Z$, $(X_s)_{s\in S}$ une famille analytique de $n$-cycles de $Z$,  $\Psi$ une famille de supports de $Z$ contenant les supports des cycles $X_s$,  $\Theta$  une famille paracompactifiante de supports de $S$  telle que la famille $\widetilde{\Psi}:=(S\times (\Phi\cap \Psi))\cap X$ soit paracompactifiante et contenue dans la famille paracompactifiante des fermés $\Theta$-propre de sorte  que le couple $(\widetilde{\Psi}, \Theta)$ soit adapté à $\pi$. Alors, il existe un morphisme d'int\'egration d'ordre sup\'erieur sur les cycles 
$${\sigma}^{q,p}_{{\Phi},X}:{\rm H}^{n+p}_{\Phi}(Z, {\mathcal L}^{n+q}_{Z})\longrightarrow 
{\rm H}^{p}_{\Theta}(S, {\mathcal L}^{q}_{S})$$
possédant les mêmes propriétés fonctorielles que ${\mathcal T}^{k}_{\pi,\mathcal L}$.\rm

\vfill\eject

\phantomsection\addcontentsline{toc}{part}{Le faisceau ${\mathcal L}^{\bullet}_{X}$.}
%
%
\svnid{$Id: S07-proof-setup.tex 269 2020-01-20 11:28:53Z kebekus $}

\section{\color{blue}{Cactérisations classiques .}}
\subversionInfo

\subsection{{Formes de première espèce.}}
\approvals{Mohamed & yes}
\par\vspace{2mm}
Si $X$ est une hypersurface à singularités isolées d'un certain 
espace numérique ${\Bbb C}^{m+1}$,
la considération de $m$-formes méromorphes n 'ayant des p\^oles qu'en 
les singularités de $X$ et
dites de {\bf première espèce} apparait déjà dans divers travaux et  en particulier dans 
[P.S] (tome 1 p.178). Rappelons qu'étant donnés
$f: ({\Bbb C}^{m+1}, 0)\rightarrow ({\Bbb C}, 0)$ un germe de 
fonction holomorphe ayant une singularité
isolée  à l'origine, un ouvert  $V:=\{z\in {\Bbb C}^{m+1} / \mid 
z\mid<\epsilon; \mid f(z)\mid<\eta\}$ et
$X$ l'ensemble d'annulation de $f$ sur $V$, on dit qu'une forme 
holomorphe $\xi$ sur 
$V-\{0\}$ est de {\bf  première espèce }\rm  s'il existe une 
désingularisation de $V$ dans laquelle $\xi$ 
se prolonge holomorphiquement.\par\noindent
Dans \cite{Gr}, Greuel a montré que toute forme holomorphe sur $V-\{0\}$ 
dont le degré est strictement inférieur
à $m-1$, est automatiquement de première espèce. Le résultat 
général  donné  par Flenner (\cite{Fl}) dit que toute $q$-forme holomorphe sur la partie 
régulière d'un espace analytique éventuellement réduit  est de  première espèce
 si  $q +1$ est  strictement inférieur à la codimension du lieu singulier.\par\noindent  On peut citer  l'article de Merle et Teissier \cite{MT} et celui de Steenbrinck et Van 
Straten \cite{SS} traitant ce sujet dans des cas plus ou moins particuliers. En degré maximum, on peut noter que ces formes méromorphes  sont caractérisées par le fait qu'elles sont de carré sommable ( ce qui implique
que leurs images directes par tout morphisme fini local de $X$ sur un 
ouvert lisse
de ${\Bbb C}^{m}$ est holomorphe) . Sur un espace ayant une singularité 
isolée, les formes de première espèce de degrés $m-1$ et $m$ jouent un 
rôle particulièrement important puisqu'elles sont attachées à  deux 
invariants fondamentaux, à
savoir, les genres géométriques et arithmétiques :
 \[p_{g}:={\rm dim}{\frac{\Gamma( V, \omega^{m}_{V})}{\Gamma( V, {\mathcal 
L}^{m}_{V})}}\]  
 \noindent qui est le nombre de conditions d'adjonctions associés à 
la singularité, et 
 \[q:={\rm dim} \frac{\Gamma( V, \omega^{m-1}_{V})}{\Gamma( V, {\mathcal 
L}^{m-1}_{V})}\]  
\noindent L'étude de ces invariants fait l'objet de plusieurs articles 
à savoir, par exemple [S.Y] auquel 
 nous renvoyons le lecteur.  Dans un contexte différent, relatif 
au degré maximal, elles interviennent aussi dans les théorèmes  
d'annulations de Grauert et Riemmenschneider \cite{GR}.
 Signalons la caractérisation suivante de \cite{SS}  dont l'analogue en 
dimension maximale se trouve énoncé dans
(\cite{MT}  p.231)\vspace{1mm}

\noindent
\subsection{{Caractérisations.}}
\approvals{Mohamed & yes}
\par\vspace{2mm}
Dans le cas d'une hypersurface, nous avons, pour les formes de degré maximal, la description suivante:
\Prop{}{}\label{Pic}(Picard [P] tome 1 p.178)\par\noindent
Soit $(X,0)$ un germe d'hypersurface réduite d' équation 
$f(x_{1},\cdots, x_{m+1}) = 0$ dans  $({\Bbb C}^{m+1},0)$ et $\Sigma$ son lieu
singulier (ou un sous ensemble rare dans $X$). Alors:\par\noindent
{\bf(i)} toute $m$ - forme méromorphe $\xi$ sur $X$ s'écrit
$$ \xi = \phi. {{dx_{2}\wedge
dx_{3}\wedge\cdots\wedge dx_{m+1}}\over{\partial{f}\over{\partial
{x_{1}}}}}$$
avec $\phi$ méromorphe sur $X$.\par\noindent
{\bf(ii)} pour que $\xi$ reste fini au voisinage de $O$, il faut que\par
a) $\phi$ soit holomorphe au voisinage de $O$.\par
b) Chaque composante irréductible $\Sigma_{i}$ de $\Sigma$, de
codimension 1, contienne un ouvert analytique dense $V_{i}$ tel que
pour tout point $x$  de $V_{i}$ ( ou représentant assez petit) on ait ${\phi.{\mathcal O}_{X,x}\subset {\mathcal C}_{x}}$ 
 avec ${\mathcal C}_{x}:= \{g\in {\mathcal O}_{X,x}; g.\widehat{{\mathcal
O}_{X,x}}\subset {\mathcal O}_{X,x}\}$ ( le {\it{conducteur }} de
l'algèbre  ${\mathcal O}_{X,x}$ dans sa normalisée $\widehat{{\mathcal
O}_{X,x}}$.
\rm\vspace{1mm} 

\noindent
On a aussi
\Prop{}{}\label{SS} (\cite{SS}, {\it proposition(1.1)},p.98)
Pour toute $q$-forme holomorphe $\xi $ sur $X-\Sigma$, les conditions 
suivantes sont équivalentes:\par
i) Pour toute application $\gamma: \Delta^{q} \rightarrow X$ de classe 
${\mathcal C}^{\infty}$ sur le
$q$- simplexe $\Delta^{q}$, l'intégrale $\int_{\gamma}\xi $ existe 
si $\gamma^{-1}(\Sigma)$ est
 d'intérieur vide.\vspace{1mm}

 \noindent ii) Il existe un morphisme de résolution des singularités $\pi 
:\tilde X\rightarrow X$ tel que
$\pi^{*}(\xi)$ se prolonge holomorphiquement à $\tilde X$ tout 
entier.\vspace{1mm}

 \noindent 
iii) L'assertion ii) est vraie  pour toute désingularisation.\vspace{1mm}

 \noindent 
iv) Il existe une modification $\pi :\tilde X\rightarrow X$, avec $\tilde 
X$ lisse, pour laquelle 
$\pi^{*}(\xi)$ se prolonge holomorphiquement à $\tilde X$  tout 
entier.\vspace{1mm}

 \noindent 
v) iv) est vraie  pour toute modification avec espace  " modifié " lisse  
\rm\vspace{1mm}

 \noindent 
\begin{rem}
{\bf(i)} On peut signaler qu'indépendamment,  Greuel (\cite{Gr}, {\it proposition(2.3)}) et Yau (\cite{Y}, {\it theorem(2.9)}) ont montré que dans le cas d'une singularité isolée d'hypersurface de dimension $m$,  toute forme génériquement holomorphe est automatiquement de première èspèce en tout degré strictement inférieur à $m-1$.  Résultat qui a été généralisé par Flenner (\cite{Fl},  {\it theorem}) en montrant que sur un espace analytique complexe (réduit) dont le lieu singulier est de codimension $r$, toute forme génériquement holomorphe est automatiquement de  première èspèce en tout degré strictement inférieur à $r-1$.\vspace{1mm}

 \noindent 
{\bf(ii)} La suite exacte courte
$$0\rightarrow{\mathcal L}^{m}_{Z}\rightarrow\omega^{m}_{Z}\rightarrow {\mathcal K}\rightarrow 0$$
où ${\mathcal K}$ est le faisceau quotient dont le support est contenu dans le lieu singulier de $Z$, nous donne la suite exacte longue
$$ 0\rightarrow{\mathcal H}om_{{\mathcal O}_{Z}}(\omega^{m}_{Z},
{\mathcal L}^{m}_{Z})\rightarrow{\mathcal H}om_{{\mathcal O}_{Z}}(\omega^{m}_{Z},
\omega^{m}_{Z})\rightarrow{\mathcal H}om_{{\mathcal O}_{Z}}(\omega^{m}_{Z},
{\mathcal K})\rightarrow \cdots$$
qui, pour  $Z$  de Gorenstein (i.e  $\omega^{m}_{Z}\simeq {\mathcal O}_{Z}$) est réduite à la suite exacte courte donnée par ces trois termes. L'identification $\omega^{m}_{Z}\simeq {\mathcal O}_{Z}$ permet de voir le faisceau ${\mathcal H}om_{{\mathcal O}_{Z}}(\omega^{m}_{Z},
{\mathcal L}^{m}_{Z})\simeq{\mathcal L}^{m}_{Z}$ comme un idéal de ${\mathcal O}_{Z}$ constitué des fonctions holomorphes sur $S$ dont la multiplication avec un générateur du faisceau dualisant $\omega^{m}_{Z}$ donne une section du faisceau ${\mathcal L}^{m}_{Z}$; c'est le {\it conducteur  d'adjonction} ${\mathcal C}$ de $Z$ (cf def. 1.6
p.238 de \cite{MT}).  Signalons que le conducteur usuel est donné par
${\mathcal H}om({\mathcal L}^{0}_{Z}, {\mathcal O}_{Z})$.\end{rem}
\vspace{1mm}

\noindent
\defn{}{} \label{def0} Un morphisme $f:X\rightarrow Y$ d'espaces complexes est dit {\emph{admissible}} si $f^{-1}({\rm Sing}(Y))$ est d'intérieur vide dans $X$.\rm

\defn{}{} \label{def0'} Un morphisme $f:X\rightarrow Y$ d'espaces complexes est dit {\emph{universellement équidimensionnel}} s'il est ouvert, surjectif à fibres de dimension ${\rm dim}(X)-{\rm dim}(Y)$.\rm

\defn{}{}\label{def1} Soit $X$ un espace analytique complexe réduit de dimension pure $m$ et de lieu singulier $\Sigma:={\rm Sing}(X)$. Si $W$ est un sous ensemble analytique de $X$ 
contenant $\Sigma$ mais ne contenant aucune composante irréductible de $X$, on dit qu'une $m$-forme holomorphe $\xi$ sur $X\setminus W$ est de type ${\rm L}^{2}$ sur $X$ si
$$\int_{U\cap(X\setminus W)}\xi\wedge\bar\xi<\infty$$
dans un voisinage $U$ de tout point de $W$. On note  ${\rm L}^{2}_{X}$ l'ensemble de ces formes.\rm\vspace{2mm}

\noindent 
Pour les degrés intermédiaires, on a naturellement 
\defn{}{}\label{def2} Avec les notations et hypothèses de la définition précédentes, on dit qu'une $k$-forme holomorphe $\xi$ sur $X\setminus W$ est de type ${\rm L}^{2}$ sur $X$ si, pour tout sous ensemble analytique $Y$ de dimension $k$ non contenu dans $\Sigma$, la restriction (générique) $\xi\mid_{Y}$ est de type ${\rm L}^{2}$ sur $Y$ au sens de la \definitionref{def1}.\rm\vspace{2mm}

\noindent 
Cette définition montre que l'on peut se contenter de tester les sous espaces $Y$ tel que ${\rm Sing}(Y)\subset W$ et éventuellement ${\rm Sing}(Y)=Y\cap {\rm Sing}(X)$.\vspace{1mm}

\noindent
Il est facile de voir  (cf par exemple \cite{GR}) que 
\Prop{}{} ${\rm L}^{2}_{X}={\mathcal L}^{m}_{X}$.
\Prop{}{} Soit $X$ un espace analytique réduit de dimension $m$. Alors, le faisceau des $k$-formes méromorphes de type ${\rm L}^{2}_{X}$ sur $X$ coïncide avec le faisceau ${\mathcal L}^{k}_{X}$.
\begin{proof} Il suffit de comparer avec la \propositionref{SS}$\,\blacksquare$
\end{proof}\vspace{1mm}

\noindent Enfin, terminons en rappelant que dans \cite{K0}, on avait montré le
 \Th{}{}\label{T0}
Soient\par\noindent
- $Z$ et $S$  deux espaces  analytiques complexes de dimension pure avec 
$Z$ dénombrable 
à l'infini et $S$ réduit.\vspace{1mm}

\noindent
- $(X_{s})_{s\in S}$ une famille analytique de $n$- cycles de $Z$ 
paramétrée par $S$, 
dont le graphe sera not\'e $X$ dans $S\times Z$.\vspace{1mm}

\noindent
- $\pi: X\rightarrow S$ est le morphisme de projection de ce graphe  sur 
$S$.\vspace{1mm}

\noindent
- $\Phi$ et  $\Psi$ deux familles de supports v\'erifiant\par
 i) $\Phi$ paracompactifiante\par
ii) $\Phi \cap \Psi$ contenue dans la famille des compacts de $Z$
\par
iii)  $\forall s\in S$, $\vert{X_{s}}\vert \in \Psi$\par
iv) $X\cap (S\times \Phi)$ contenu dans la famille des ferm\'es $S$- 
propres.\vspace{1mm}

\noindent  
 Alors il existe
un morphisme d'int\'egration d'ordre sup\'erieur sur les cycles 

$$\tilde{\sigma}^{q,0}_{{\Phi},X}:{\rm H}^{n}_{\Phi}(Z, {\mathcal 
L}^{n+q}_{Z})\longrightarrow 
{\rm H}^{0}(S, {\mathcal L}^{q}_{S})$$
\noindent v\'erifiant les propri\'et\'es suivantes : \vspace{1mm}

\noindent
a) {\bf Compatibilit\'e avec les inclusions ouvertes sur} 
$Z$:\vspace{1mm}

\noindent
Soit $U$ un ouvert de $Z$ et $S'$ un ouvert de $S$ form\'e des points 
pour 
lesquels les cycles associ\'es sont dans $U$. Si  $j:U\rightarrow Z$ est 
l'injection naturelle, on a
 la commutativit\'e du diagramme suivant:
$$\xymatrix{ {\rm H}^{n}_{\Phi}(Z, {\mathcal L}^{n+q}_{Z})\ar 
[d]_{\tilde{\sigma}^{q,0}_{{\Phi},X}} 
& \ar[l] {\rm H}^{n}_{{\Phi}\mid U}(U, {\mathcal L}^{n+q}_{U})\ar 
[d]^{{\tilde{\sigma}^{q,0}_{{\Phi},X}}\mid U}\cr
 {\rm H}^{0}(S, {\mathcal L}^{q}_{S} ) \ar[r]& {\rm H}^{0}(S', {\mathcal L}^{q}_{S'} )}$$
dans lequel les fl\`eches horizontales sup\'erieure et inf\'erieure 
d\'esignent respectivement
 le prolongement  par zero des sections \`a supports  et la restriction  
usuelle.\vspace{1mm}

\noindent   
b) {\bf Compatibilit\'e avec les images directes de cycles}\vspace{1mm}

\noindent
Soit $f: Z\rightarrow T$ est un morphisme propre et surjectif d'espaces 
analytiques complexes 
 r\'eduits et $\tilde\Phi$ une famille  paracompactifiante de ferm\'es de 
$T$ v\'erifiant la condition 
iv) ci-dessus  pour la famille analytique $(f_{*}X_{s})_{s\in S}$ dont le 
graphe sera 
not\'e $\hat{X}$ et telle que $f^{-1}\tilde\Phi$ soit contenue dans 
$\Phi$, alors on a le diagramme commutatif :
  $$\xymatrix{ {\rm H}^{n}_{\tilde\Phi}(T, {\mathcal L}^{n+q}_{T})
\ar[dr]_{{\tilde{\sigma}}^{q,0}_{{\tilde\Phi},\hat{X}}} 
\ar [rr]^{f^{*}} && 
 {\rm H}^{n}_{\Phi}(Z, {\mathcal L}^{n+q}_{Z})\ar 
[dl]^{{\tilde{\sigma}}^{q,0}_{{\Phi},X}}\cr
 &{\rm H}^{0}(S, {\mathcal L}^{q}_{S} )} $$
 \noindent 
c) {\bf Compatibilit\'e avec les changements de bases sur $S$} 
:\vspace{1mm}

\noindent
   Soit $\phi: \tilde{S}\rightarrow S$ un morphisme d'espaces analytiques 
complexes r\'eduits. 
Si $\tilde{X}$ d\'esigne le graphe de la famille
$(X_{\tilde s})_{\tilde {s}\in\tilde { S}}$ obtenue par changement de 
base sur $S$, on a le diagramme commutatif suivant:
$$\xymatrix{ {\rm H}^{n}_{\Phi}(Z, {\mathcal L}^{n+q}_{Z})\ar 
[d]^{_{\tilde{\sigma}^{q,0}_{{\Phi},X}}}  
\ar [dr]^{_{\tilde{\sigma}^{q,0}_{{\Phi},\tilde{X}}}} \cr 
 {\rm H}^{0}(S, {\mathcal L}^{q}_{S} )\ar [r]_{\phi^{* }}& {\rm H}^{0}(\tilde{S}, {\mathcal 
L}^{q}_{\tilde{S}} )} $$
\noindent
d) {\bf Faisceautisation:}\vspace{1mm}

\noindent
$\pi$ induit un morphisme de faisceau de ${\mathcal O}_{S}$ modules que l'on 
appelle trace 
{\it faisceautique}
$$ {\rm I}\!{\rm R}^{n}\pi_{!}{\mathcal L}^{n+q}_{X}\rightarrow {\mathcal L}^{q}_{S}$$
\noindent fonctorielle en $S$ et en $X$, qui d\'etermine et qui est 
enti\`erement d\'etermin\'e, 
dans ${\bf D}({\mathcal O}_{S})$ ( la cat\'egorie d\'eriv\'ee des ${\mathcal 
O}_{S}$-modules) par la fl\`eche
$$ {\rm I}\!{\rm R}\pi_{!}{\mathcal L}^{n+q}_{X}[n]\rightarrow {\mathcal L}^{q}_{S}. $$
\noindent 
e) {\bf Compatibilit\'e avec " les cas usuels "}\par
i) Pour $Z$ et $S$ lisses, sa construction coincide avec celle faite en 
termes d'image directe au sens des courants.
\par 
ii)  Pour $q=0$, le lien avec le morphisme $\rho^{0}_{{\Phi},X}$ 
construit dans \cite{B.V} est donn\'e par 
 le diagramme commutatif

$$\xymatrix{ {\rm H}^{n}_{\Phi}(Z, \Omega^{n}_{Z})\ar[d]^{\rho^{0}_{\Phi,X}} 
 \ar[r]&H^{n}_{\Phi}(Z, {\mathcal L}^{n}_{Z})\ar[d]_{\tilde{\sigma}^{0,0}_{{\Phi},X}}\\
 {\rm H}^{0}(S, {\mathcal O}_{S} )\ar[r]&{\rm H}^{0}(S, {\mathcal L}^{0}_{S} ) }$$
 \noindent 
dans lequel les fl\`eches horizontales sup\'erieure et inf\'erieure 
d\'ecoulent respectivement
 des morphismes canoniques $\Omega^{n}_{Z}\rightarrow {\mathcal L}^{n}_{Z}$ 
et 
${\mathcal O}_{S}\rightarrow {\mathcal L}^{0}_{S}$.\rm
\vspace{2mm}

\noindent 
On avait fait remarquer que, pour $q=1$, le morphisme d'int\'egration est à valeurs dans un sous faisceau de ${\mathcal L}^{1}_{S}$, 
$$ {\rm H}^{n}_{\Phi}(Z, \Omega^{n+1}_{Z})\rightarrow {\rm H}^{0}(S, \overline{\Omega}^{1}_{S})$$ 
 $ \overline{\Omega}^{1}_{S}$ pouvant être décrit comme étant le faisceau des formes m\'eromorphes
 v\'erifiant une relation de d\'ependance int\'egrale sur l'algèbre symétrique des formes holomorphes, ou des formes m\'eromorphes, qui sur le tangeant de Zariski, correspondent aux fonctions m\'eromorphes localement born\'ees (cf \cite{K3} pour une discussion approfondie et détaillée). 
 \vspace{1mm}

 \noindent 
\section{\color{blue}{Quelques résultats basiques.}}
\subsection{Propriétés fondamentales.}
Dans ce paragraphe, nous donnons quelques résultats qui, pour la plupart, se réfèrent à \cite{K3} où la fonctorialité du faisceau $\omega^{\bullet}_{X}$ et ses images directes supérieures sont étudiées en détails. 
\Prop{}{}\label{P1} Soit $X$ un espace analytique complexe réduit et de dimension pure $m$. Alors, pour tout entier $k\geq 0$, on a:\vspace{1mm}

\noindent
{\bf(i)}  le faisceau ${\mathcal L}^{k}_{X}$  est un sous faisceau   de $\omega^{k}_{X}$  (en particulier, il est  sans  ${\mathcal O}_{X}$- torsion et  coincide avec le faisceau des formes holomorphes usuelles aux points réguliers de $X$). De plus, ces faisceaux sont  stables par différentiation extérieure et  munis d'un cup-produit interne. \vspace{1mm}

\noindent 
{\bf(ii)}  tout morphisme $f:X\rightarrow Y$ morphisme fini et surjectif d'espaces analytiques complexes réduits de dimension pure $m$ induit, pour tout entier $k\leq m$, un morphisme {\it{trace}}
$$f_{*}{\mathcal L}^{k}_{X}\rightarrow {\mathcal L}^{k}_{Y}$$
{\bf(iii)} Toute modification propre $f:X\rightarrow Y$ induit un isomorphisme 
$$f_{*}{\mathcal L}^{k}_{X}\simeq {\mathcal L}^{k}_{Y}$$
{\bf(iv)} si ${\mathcal F}$ est un sous faisceau cohérent de $\omega^{k}_{X}$ coincidant avec le faisceau des formes usuelles sur la partie régulière et pour lequel il existe {\bf un} morphisme propre $f: Z\rightarrow X$ avec $Z$ lisse tel que l'image réciproque par $f$ de toute section de ${\mathcal F}$ se prolonge holomorphiquement sur $Z$, alors ${\mathcal F}$ est un sous faisceau de ${\mathcal L}^{k}_{X}$.\vspace{1mm}

\noindent 
\rm
\begin{proof}\vspace{1mm}

\noindent
{\bf(i)}  Comme, par définition, ${\mathcal L}^{k}_{X}=\pi_{*}(\Omega^{k}_{\tilde X})$ pour toute désingularisation $\pi:\tilde{X}\rightarrow X$, il est cohérent puisqu'image directe propre d'un faisceau cohérent (d'après le théorème de cohérence de Grauert). Il est aussi, de ce fait, sans torsion  puisque $\Omega^{k}_{\tilde X}$ l'est. De plus,  pour tout ouvert $U$ de $X$, on a   ${\mathcal L}^{k}_{X}|_{U}= {\mathcal L}^{k}_{U}$ et   coincide naturellement avec  le faisceau des formes holomorphes usuelles aux points réguliers de $X$ puisque $\pi$ est un isomorphisme sur cet ouvert.  Comme l'image directe au sens des courants donne un morphisme injectif $\pi_{*}\Omega^{k}_{\tilde X}\rightarrow \omega^{k}_{X}$,  il est de profondeur au moins $1$;  en particulier, pour  tout sous espace $T$  rare ou d'intérieur vide dans $X$, on a l'annulation des faisceaux de cohomologie ${\mathcal H}^{0}_{T}({\mathcal L}^{k}_{X})=0$.\vspace{1mm}

\noindent
Pour définir le cup produit, il suffit d'écrire
$$ {\mathcal L}^{k}_{X}\otimes{\mathcal L}^{k'}_{X}= \pi_{*}\Omega^{k}_{\tilde X} \otimes \pi_{*}\Omega^{k'}_{\tilde X}  \rightarrow  \pi_{*}(\Omega^{k}_{\tilde X} \otimes \Omega^{k'}_{\tilde X} )\rightarrow \pi_{*}( \Omega^{k+k'}_{\tilde X})={\mathcal L}^{k+k'}_{X}$$
{\bf(ii)} Cette trace se construit grâce à la trace déjà vue $\pi_{*}\omega^{k}_{X}\rightarrow\omega^{k}_{Y}$ induite par l'image directe au sens des courants. Ainsi, par composition avec le morphisme canonique ${\mathcal L}^{k}_{X}\rightarrow\omega^{k}_{X}$ et pour $Y$ lisse, on a un morphisme  $\pi_{*}{\mathcal L}^{k}_{X}\rightarrow\Omega^{k}_{Y}$. Le cas général s'en déduit immédiatement en considérant une désingularisation de $Y$ et le diagramme de changement de base associé
$$\xymatrix{\tilde{X}\ar[r]^{\tilde\phi}\ar[d]_{\tilde{\pi}}&{X}\ar[d]^{\pi}\\ 
 {\tilde Y}\ar[r]_{\phi}& Y}$$
 Alors, comme $\pi_{*}{\mathcal L}^{k}_{Z}=\pi_{*}{\tilde\phi}_{*}{\mathcal L}^{k}_{\tilde Z}={\phi}_{*}{\tilde \pi}_{*}{\mathcal L}^{k}_{\tilde Z}$, il nous suffit d'appliquer le foncteur $\phi_{*}$ au morphisme trace ${\tilde \pi}_{*}{\mathcal L}^{k}_{\tilde X}\rightarrow \Omega^{k}_{\tilde Y}$. 
 Rappelons que le morphisme trace commute à tout changement de base c'est-à-dire 
$$\phi^{*}{\rm Tr}_{\pi}(\xi)={\rm Tr}_{\tilde\pi}({\tilde\phi}\xi)$$
 {\bf(iii)} Cela découle immédiatement de la définition.\vspace{1mm}

\noindent
{\bf(iv)} Soient ${\mathcal F}$ un sous faisceau cohérent de $\omega^{k}_{X}$ coincidant avec le faisceau des formes holomorphes  usuelles sur la partie régulière et pour lequel il existe {\bf un} morphisme propre $f: Z\rightarrow X$ avec $Z$ lisse tel que l'image réciproque par $f$ de toute section de ${\mathcal F}$ se prolonge holomorphiquement sur $X$. Montrons que c'est nécessairement un sous faisceau du faisceau ${\mathcal L}^{k}_{X}$. \vspace{1mm}

\noindent
Commençons par remarquer que dans le cas où $f$ est une modification propre, le résultat est clair. De même, si $f$ est fini et surjectif, le changement de base donnée par la normalisation de $X$, nous ramène au cas d'un  revêtement ramifié de degré $r$ duquel résulte l'assertion. Précisément, relativement au diagramme de changement de base
$$\xymatrix{\tilde{Z}\ar[r]^{\tilde \phi}\ar[d]_{\tilde{f}}&{Z}\ar[d]^{f}\\ 
 {\tilde X}\ar[r]_{\phi}& X}$$
 et au morphisme trace ${\tilde f}_{*}\Omega^{k}_{\tilde Z}\rightarrow{\mathcal L}^{k}_{X}$,    on voit que la forme méromorphe ${\tilde f}_{*}({\tilde f}^{*}(\phi^{*}\xi)$ définit une section du faisceau ${\mathcal L}^{k}_{\tilde X}$ c'est-à-dire $\phi^{*}\xi$ section de ${\mathcal L}^{k}_{\tilde X}$ puisque ${\tilde f}_{*}{\tilde f}^{*}=r.{\rm Id}$. On définit ainsi un morphisme non trivial $\phi^{*}{\mathcal F}\rightarrow{\mathcal L}^{k}_{\tilde X}$. D'où, en composant avec la flèche naturelle ${\mathcal F}\rightarrow\phi_{*}\phi^{*}{\mathcal F}$, le morphisme
 $${\mathcal F}\rightarrow{\mathcal L}^{k}_{ X}$$
qui est nécessairement injectif puisque c'est un isomorphisme générique (les deux faisceaux coincidant avec le faisceau des formes holomorphes sur la partie lisse) et  ${\mathcal F}$ est sans torsion car sous faisceau du faisceau sans torsion  $\omega^{k}_{X}$.\vspace{1mm}

\noindent 
Remarquons, au passage, que $\phi$ étant le morphisme de normalisation, on a 
 $$\phi_{*}\phi^{*}({\mathcal L}^{k}_{X})={\mathcal L}^{k}_{X}\otimes \phi_{*}{\mathcal O}_{\tilde X}={\mathcal L}^{k}_{X}\otimes {\mathcal L}^{0}_{X}={\mathcal L}^{k}_{X}$$ 
 Dans le cas général, on peut utiliser\footnote{On peut prendre pour $\phi$ une désingularisation de $X$ et  raisonner en termes de courants en utilisant la formule de projection ${{\tilde f}_{*}}([\tilde{Z}]\wedge {\tilde f}^{*}(\phi^{*}\xi))=[\tilde{X}]\wedge {\phi}^{*}(\xi)$ et voir ${\phi}^{*}(\xi)$ comme  un courant de type $(k,0)$ et $\bar\partial$-fermé donc définissant, en vertu du lemme de Dolbeault-Grothendieck,  une $k$-forme holomorphe sur $\tilde X$ et on termine comme dans le cas fini.} le théorème d'applatissement de Hironaka \cite{H1} pour avoir un diagramme similaire au précédent avec $f$ plat, $\phi$ et $\tilde\phi$ des mordifications et $\tilde Z$ la réunion des composantes irréductibles du produit fibré ${\tilde X}\times_{X} Z$ qui se surjectent sur $\tilde X$. On peut, alors, tariter le problème localement et nous ramener essentiellement au cas fini.\vspace{1mm}

 \noindent 
  Supposons, donc, $f: Z\rightarrow X$ soit un morphisme plat d'espaces complexes réduits  et considérons  une factorisation locale $\xymatrix{ Z\ar@/^1pc/[rr]^{f}\ar[r]_{h}& X\times U\ar[r]_{q}&X}$ de $f$. Comme  $f^{*}(\xi)=h^{*}(q^{*}(\xi))$ est holomorphe, sa trace ${\rm Tr}_{h}(h^{*}(q^{*}(\xi))$ définit une section du faisceau ${\mathcal L}^{k}_{X\times U}$; ce qui signifie exactement  que $q^{*}(\xi)$ est une section de ${\mathcal L}^{k}_{X\times U}$  d'être une section du faisceau ${\mathcal L}^{k}_{Z\times U}$. Mais comme $q$ n'est rien d'autre que la projection canonique, cela revient exactement à dire  que $\xi$ est une section du faisceau ${\mathcal L}^{k}_{X}$.
Traduisant ceci dans notre contexte, cela nous dit que pour toute section $\xi$ de ${\mathcal F}$,  $\phi^{*}\xi$ est une section de ${\mathcal L}^{k}_{\tilde X}$. D'où le morphisme $\phi^{*}{\mathcal F}\rightarrow {\mathcal L}^{k}_{\tilde X} $ et, donc, ${\mathcal F}\rightarrow {\mathcal L}^{k}_{X} $ nécessairement injectif pour les raisons évoquées dans le cas fini ci-dessus$\,\blacksquare$
\end{proof}
\Prop{}{}\label{P2} Soit $X$ un espaces analytique complexe réduit et soit $\xi$ une section du faisceau $\omega^{k}_{X}$. Alors, $\xi$ est une section du faisceau ${\mathcal L}^{k}_{X}$ si et seulement si il existe une variété complexe lisse $Y$ et un morphisme $\pi:Y\rightarrow X$  propre et surjectif tel que $\pi^{*}(\xi)$ génériquement holomorphe sur $Y$ se prolonge holomorphiquement sur $Y$.\rm
\begin{proof}
$\bullet$ Si $\xi$ est une section du faisceau ${\mathcal L}^{k}_{X}$, alors, par définition, il existe une résolution des singularités de $X$ donc un morphisme propre et surjectif $\pi:\widetilde{X}\rightarrow X$ tel que $\pi^{*}(\xi)$ se prolonge holomorphiquement sur $\widetilde{X}$.\vspace{1mm}

\noindent
$\bullet$ Pour la réciproque,  le diagramme  commutatif 
$$\xymatrix{Z\ar[d]_{\tilde \pi}\ar[rd]_{\psi}\ar[r]^{\theta}&Y\ar[d]^{\pi}\\
\widetilde{X}\ar[r]_{\nu}&X}$$
dans lequel $\widetilde X$ est une résolution de $X$ et $Z$ une résolution du produit fibré réduit $\widetilde{X}\times_{X} Y$, nous ramène à traiter essentiellement  la situation donnée par le diagramme commutatif
$$\xymatrix{Y\ar[d]_{\psi}\ar[rd]_{\pi}\\
\widetilde{X}\ar[r]_{\nu}&X}$$
Soit $\xi$ une forme méromorphe sur $X$ (à pôles dans le lieu singulier de $X$). Alors, $\nu^{*}(\xi)$ est une forme méromorphe à pôles dans le diviseur exceptionnel dont l'image réciproque $\psi^{*}(\nu^{*}(\xi))$  est holomorphe sur $Y$. Mais ceci n'est possible que si et seulement si $\nu^{*}(\xi)$ est déjà holomorphe sur $\widetilde{X}$ comment on s'en convainc aisément en utilisant le lemme de Dolbeault-Grothendieck. En effet, $\pi^{*}(\xi)$ étant holomorphe définit un courant $\overline\partial$-fermé dont l'image directe $\psi_{*}(\pi^{*}(\xi))$ est encore un courant $\overline\partial$-fermé. Mais, $\widetilde X$ étant lisse, ce courant correspond exactement à une forme holomorphe. On termine, alors,  en constatant que
$$\psi_{*}(\pi^{*}(\xi))=\psi_{*}\psi^{*}(\nu^{*}\xi)=\nu^{*}\xi$$
puisque $\widetilde X$ est lisse (donc normal!)\footnote{ On peut supposer $X$ normal et utiliser la factorisation de Stein de $\pi$ donnée par $\xymatrix{Y\ar@/^1pc/[rr]^{\pi}\ar[r]_{\alpha}&Y'\ar[r]_{\beta}&X}$  avec $\alpha$ propre à fibres connexes (i.e ${\alpha_{*}}{\mathcal O}_{Y}={\mathcal O}_{Y'}$) et $\beta$ fini, surjectif et ouvert puisque $X$ est normal. alors, $\pi_{*}\pi^{*}(\xi)=\beta_{*}(\alpha_{*}\alpha^{*}(\beta^{*}(\xi)))=\beta_{*}\beta^{*}(\xi)$ supposé holomorphe par Dolbeault-Grothendieck. Mais comme $\beta$ est un revêtement ramifié d'un certain degré $k$, le courant $\beta_{*}\beta^{*}(\xi)$ coïncide avec la trace ${\mathcal T}r_{\beta}(\beta^{*}\xi)=k.\xi$ qui montre bien que $\xi$ est holomorphe.}\vspace{1mm}

\noindent On en déduit, alors,  que   $\xi$ est une section du faisceau ${\mathcal L}^{k}_{X}\,\blacksquare$ \end{proof}
\vspace{1mm}

\noindent 
\cor{}{}\label{cor2} Soient $k$ un entier naturel non nul et $X$ un espace analytique complexe réduit.  Alors,  le faisceau $\omega^{k}_{X}$ est stable par image réciproque de morphismes d'espaces complexes si et seulement si $\omega^{k}_{X}={\mathcal L}^{k}_{X}$.\rm\par\noindent
\begin{proof} Supposons que ce soit le cas et considérons, à cet effet,  le morphisme de désingularisation $\pi:\tilde{X}\rightarrow X$.  Alors, par hypothèse, on a un morphisme  $\omega^{k}_{X}\rightarrow \pi_{*}\omega^{k}_{\tilde X}$. Comme $\pi_{*}\omega^{k}_{\tilde X}=\pi_{*}\Omega^{k}_{\tilde X}={\mathcal L}^{k}_{X}$,    $\omega^{k}_{X}\rightarrow\pi_{*}\Omega^{k}_{\tilde X}={\mathcal L}^{k}_{X}$, on a, donc, un morphisme $\omega^{k}_{X}\rightarrow \rightarrow {\mathcal L}^{k}_{X}$ qui est nécessairement injectif puisque   ces faisceaux coincident génériquement sur $X$ et sont sans torsion.  Or, on sait déjà que ${\mathcal L}^{k}_{X}$ est un sous faisceau de $\omega^{k}_{X}$. D'où la conclusion$\,\,\blacksquare$
\end{proof}
\begin{rem}
La condition de surjectivité est essentielle pour éviter le cas d'un plongement pour lequel le critère n'est pas valable comme le montre l'exemple simple suivant
de $X:=\{(x,y,z)\in {\Bbb C}^{3}: xy^{2}=z^{2}\}$ et de la forme $\displaystyle{\frac{dz}{y}}$ qui est une section du faisceau $\omega^{1}_{X}$ sans l'être pour ${\mathcal L}^{1}_{X}$ et qui, par restriction à la courbe ($x=t^{2}, y=t, z=t^{2}$), donne une forme holomorphe!
\end{rem}
Signalons le résulat relativement récent de \cite{KeSc} ({\bf{theorem 1.1}}) que l'on peut reécrire sous la forme
\Th{}{} \label{KS} Soit $X$ un  espace analytique complexe réduit de dimension pure $m$. Alors, une section $\xi$ du faisceau $\omega^{k}_{X}$ est une section du faisceau  ${\mathcal L}^{k}_{X}$ si et seulement si $\xi$ (resp. $d\xi$) est une section de ${\mathcal H}om(\Omega^{m-k}_{X}, {\mathcal L}^{m}_{X})$ (resp. de ${\mathcal H}om(\Omega^{m-k-1}_{X}, {\mathcal L}^{m}_{X})$).\rm

\phantomsection\addcontentsline{toc}{part}{  Le \theoremref{T1} pour certains morphismes.}
%
%
\svnid{$Id: S07-proof-setup.tex 269 2020-01-20 11:28:53Z kebekus $}

\section{\color{blue}{Image réciproque pour  certaines  classes de morphismes.}}
\subversionInfo
\approvals{Mohamed & yes}
\par\vspace{2mm}

\noindent 
Dans ce paragraphe, nous mettons en évidence quelques situations pour lesquelles la construction d'une image réciproque ne nécessite que très peu d'efforts.  L'aspect fonctoriel des constructions a été volontairement survolé sachant qu'il sera développé dans le cas général.
\Th{}{}\label{T3} Soit $\pi:X\rightarrow S$ un morphisme 
 d'espaces analytiques réduits de dimension pure $m$ et $r$ respectivement. Alors, si\vspace{1mm}
 
 \noindent
 {\bf(i)} $\pi$ est surjectif avec  l'une des propriétes: propre, universellement $(m-r)$-équidimensionnel ou $(m-r)$-équidimensionnel ou  \vspace{1mm}

\noindent 
 {\bf(ii)} $\pi$ est admissible (i.e ${\pi}^{-1}({\rm Sing}(S))$ d'intérieur vide dans $X$),\vspace{1mm}

 \noindent
 il existe un unique   morphisme de faisceaux cohérents
$\displaystyle{{\bar\pi}^{*}:{\mathcal L}^{\bullet}_{S}\rightarrow \pi_{*}{\mathcal L}^{\bullet}_{X} }$ prolongeant naturellement l'image réciproque des formes holomorphes usuelles, compatible avec la composition des morphismes de même nature que dans {\bf(i)} ou {\bf(ii)} c'est-à-dire si
$\xymatrix{X\ar@/^1pc/[rr]^{h}\ar[r]_{f}&Y\ar[r]_{g}&Z}$ est la composée de morphismes d'espaces complexes réduits avec $f$ et $g$
surjectifs propres ou $m-r$-équidimensionnel ou $m-r$-équidimensionnel ou bien satisfaisant les relations d'incidence 
$g^{-1}({\rm Sing}(Z))$ (resp. $f^{-1}({\rm Sing}(Y)\cup h^{-1}({\rm Sing}(Z))$) soit d'intérieur vide dans $Y$ (resp. $X$), on a
$\displaystyle{{\bar h}^{*}={\bar f}^{*}\circ {\bar g}^{*}}$.
\rm \vspace{1mm}

\noindent 
\begin{proof}
{\bf{(i)}} \vspace{1mm}
\indent $\bullet$ {\bf{Cas propre et surjectif:}} Il peut être vu comme un corollaire de la \propositionref{P2} puisque pour une résolution des singularités $\theta:\tilde{X}\rightarrow X$, la composée $\psi:=\pi\circ \theta$ est un morphisme propre surjectif d'une variété lisse sur $S$. D'après cette proposition, l'image réciproque $\psi^{*}(\xi)$ de toute section $\xi$ du faisceau ${\mathcal L}^{k}_{S}$  se prolonge holomorphiquement sur $\tilde{X}$ c'est-à-dire  $\theta^{*}(\pi^{*}(\xi))$ holomorphe. Mais cela signifie exactement que $\pi^{*}(\xi)$ est une section du faisceau ${\mathcal L}^{k}_{X}$. On définit ainsi le morphisme de ${\mathcal O}_{S}$-modules cohérents
$${\bar\pi}^{*}:{\mathcal L}^{\bullet}_{S}\rightarrow \pi_{*}{\mathcal L}^{\bullet}_{X}$$
défini par ${\bar\pi}^{*}(\xi):=\theta_{*}(\psi^{*}(\xi))$ et 
dépendant à priori de $\theta$. Mais, on en vérifie facilement l'indépendance. En effet,  sachant que deux résolutions peuvent être coiffées par une troisième, on se ramène à des diagrammes commutatifs du type
$$\xymatrix{&X_{3}\ar[dd]_{\phi"}\ar[ld]_{\phi}\ar[rd]^{\phi'}&\\
X_{1}\ar[rd]_{\theta}&&X_{2}\ar[ld]^{\theta'}\\
&X\ar[d]_{\pi}&\\
&S&}$$
dont il est facile de déduire la conclusion grâce au lemme de Dolbeault-Grothendieck  puisque 
$$\phi_{*}\phi^{*}={\rm Id},\,\, \phi'_{*}\phi'^{*}={\rm Id}$$
Comme on peut le remarquer, c'est essentiellement les mêmes arguments que ceux évoqués pour l'indépendance du faisceau ${\mathcal L}^{k}_{X}$ vis-à-vis de la résolution choisie.\vspace{1mm}

\noindent 
Pour établir la {\emph{compatibilité avec la composition}} des morphismes propres surjectifs, on considère la composée donnée par $\xymatrix{ X\ar@/^1pc/[rr]^{h}\ar[r]_{f}& Y\ar[r]_{g}&S}$ et un  diagramme du type
$$
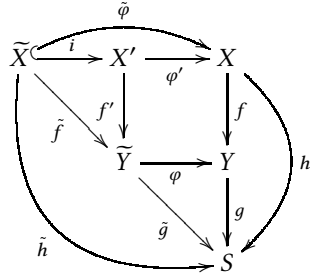
\xymatrix{\widetilde{X}\ar@{^{(}->}[r]^{i}\ar@/_3pc/[rrdd]_{\tilde h}\ar@/^1pc/[rr]^{\tilde\phi}\ar[rd]_{\tilde{f}}&X'\ar[r]_{\phi'}\ar[d]_{f'}&{X}\ar[d]^{f} \ar@/^2pc/[dd]^{h}\\ 
 &{\widetilde Y}\ar[rd]_{\tilde g}\ar[r]_{\phi}& Y\ar[d]^{g}\\
 &&S}$$
 où $X'$ (resp. $\widetilde X$) désigne la préimage totale (resp. stricte) pour la désingularisation $\phi$ de $Y$, et  de voir que la relation $(\overline{f\circ g})^{*}=\overline{g}^{*}\circ\overline{f}^{*}$ se déduit de la comparaison des courants
 $$\overline{h}^{*}(\xi):={\tilde\phi}_{*}{\tilde\phi}^{*}(h^{*}(\xi))={\tilde\phi}_{*}{\tilde\phi}^{*}(f^{*}(g^{*}(\xi)))={\tilde\phi}_{*}({\tilde f}^{*}(\phi^{*}(g^{*}(\xi))))$$
$$\overline{f}^{*}(\overline{g}^{*}(\xi)):={\tilde\phi}_{*}{\tilde\phi}^{*}(f^{*}({\phi}_{*}{\phi}^{*}g(\xi)))={\tilde\phi}_{*}({\tilde f}^{*}(\phi^{*}({\phi}_{*}{\phi}^{*}g^{*}(\xi))))$$
 puisque, en tant que courants, on a ${\phi}_{*}{\phi}^{*}g^{*}(\xi))=g^{*}(\xi)$ car, pour toute forme ${\mathcal C}^{\infty}$ de type $(m-k, m)$ sur $Y$, 
 $$\langle{{\phi}_{*}{\phi}^{*}g^{*}(\xi), \psi}\rangle=\langle{{\phi}^{*}g^{*}(\xi), \phi^{*}\psi}\rangle=\int_{{\widetilde{Y}}\setminus {\rm E}}{\phi}^{*}g^{*}(\xi))\wedge\phi^{*}\psi=\int_{{\reg}(Y)}g^{*}(\xi)\wedge \psi=\langle{g^{*}(\xi), \psi}\rangle$$
$\bullet$ {\bf{Cas d'un morphisme à fibres de dimension constante ouvert ou équidimensionnel.}}  Comme les faisceaux considérés sont invariants par modifications propres et que l'ouverture et la constance de la dimension des fibres sont des notions  préservées dans tout changement de base, on peut supposer la base normale. En effet, si  $\nu:\widetilde{S}\rightarrow S$ est la normalisation de $S$ et 
$$\xymatrix{\widetilde{X}\ar[r]^{\tilde\nu}\ar[d]_{\tilde{\pi}}&{X}\ar[d]^{\pi}\\ 
 {\widetilde S}\ar[r]_{\nu}& S}$$
 le diagramme  induit par changement de base ($\nu$ et $\tilde\nu$ étant des modifiactaions propres), on voit que l'existence d'une image réciproque pour $\tilde{\pi}$ entraine automatiquement celle pour $\pi$ puisqu'il suffit  d'appliquer le foncteur $\nu_{*}$  au  morphisme
 ${\mathcal L}^{k}_{\widetilde S}\rightarrow \tilde\pi_{*}{\mathcal L}^{k}_{\widetilde X}$ pour avoir le morphisme $\nu_{*}{\mathcal L}^{k}_{\widetilde S}\rightarrow \nu_{*}\tilde\pi_{*}{\mathcal L}^{k}_{\widetilde X}$ qui  s'écrit aussi ${\mathcal L}^{k}_{S}\rightarrow \pi_{*}{\mathcal L}^{k}_{X}$.\vspace{1mm}

 \noindent 
On peut, alors,  utiliser le {\bf{théorème 1}} de \cite{K3} et la caractérisation  de \cite{KeSc} disant que, pour un espace complexe réduit $Y$ de dimension pure $m$, 
$$\xi\in {\mathcal L}^{k}_{Y}\,\Longleftrightarrow\,\xi\wedge \alpha\,\,{\rm et}\, d\xi\wedge \beta\in {\mathcal L}^{m}_{Y},\,\forall\,\alpha\in\Omega^{m-k}_{Y},\,\,\forall\,\beta \in\Omega^{m-k-1}_{Y}$$
se traduisant, vue l'injection naturelle $${\mathcal L}^{k}_{Y}\into{\mathcal H}om({\Omega}^{m-k}_{Y}, {\mathcal L}^{m}_{Y})$$ 
par\vspace{1mm}

\noindent 
\emph{une section $\xi$ de ${\mathcal H}om({\Omega}^{m-k}_{Y}, {\mathcal L}^{m}_{Y})$ est une section du faisceau ${\mathcal L}^{k}_{Y}$ si et seulement si  $d\xi$ est une section de ${\mathcal H}om({\Omega}^{m-k-1}_{Y}, {\mathcal L}^{m}_{Y})$ (stabilité par différentiation extérieure).} \vspace{2mm}

\noindent \rm 
On notera $\widetilde{\mathcal L}^{k}_{Y}:= {\mathcal H}om({\mathcal L}^{m-k}_{Y}, {\mathcal L}^{m}_{Y})$ le sous faisceau naturel de ${\mathcal H}om({\Omega}^{m-k}_{Y}, {\mathcal L}^{m}_{Y})$ . \vspace{1mm}

\noindent
Alors, d'après le  {\bf{théorème 1}} de \cite{K3}, on a un morphisme d'image réciproque défini en tout degré $k$:
$$\pi^{*}\widetilde{\mathcal L}^{k}_{S}\rightarrow \widetilde{\mathcal L}^{k}_{X} \,{\rm ou}\,\,\widetilde{\mathcal L}^{k}_{S}\rightarrow \pi_{*}\widetilde{\mathcal L}^{k}_{X}$$
Alors, comme toute section $\xi$ de ${\mathcal L}^{k}_{S}$ est aussi une section $\widetilde{\mathcal L}^{k}_{S}$ (puisque   ${\mathcal L}^{k}_{S}\subset\widetilde{\mathcal L}^{k}_{S} $), $\pi^{*}(\xi)$ définit une section de $\widetilde{\mathcal L}^{k}_{X} $ et, donc, de  ${\mathcal H}om({\Omega}^{m-k}_{X}, {\mathcal L}^{m}_{X})$. \vspace{1mm}

\noindent Par ailleurs, 
$$\xi\in {\mathcal L}^{k}_{S}\,\Longrightarrow\,d\xi \in {\mathcal L}^{k+1}_{S}$$
et, par conséquent, $$\pi^{*}(d\xi):=d\pi^{*}(\xi)\in\widetilde{\mathcal L}^{k+1}_{X}\subset{\mathcal H}om({\Omega}^{m-k-1}_{X}, {\mathcal L}^{m}_{X})$$  
Cela prouve, d'après le critère de \cite{KeSc}, que $\pi^{*}(\xi)$ définit une section du faisceau $ {\mathcal L}^{k}_{X}$ et, par suite, l'existence de diagrammes commutatifs
$$\xymatrix{\widetilde{\Omega}^{k}_{S}\ar[d]\ar@{^{(}->}[r]&{\mathcal L}^{k}_{S}\ar@{^{(}->}[r]\ar[d]&\widetilde{\mathcal L}^{k}_{S} \ar[d]\\
\pi_{*}\widetilde{\Omega}^{k}_{X}\ar@{^{(}->}[r]&\pi_{*}{\mathcal L}^{k}_{X}\ar@{^{(}->}[r]&\pi_{*}\widetilde{\mathcal L}^{k}_{X} }$$ 
$\widetilde{\Omega}^{k}_{S}$ désignant le faisceau $\Omega^{k}_{S}$ modulo torsion (idem sur $X$); assurant la compatibilité avec l'image réciproque des formes holomorphes usuelles.
\vspace{1mm}

\noindent
le cas d'un morphisme propre surjectif peut aussi se déduire du cas géométriquement plat  en vertu du théorème d'applatissement algébrique\footnote{ Dans \cite{B2}, on trouve une version  géométrique traitant le cas d'un morphisme propre génériquement à fibres de dimension constante et génériquement réduites et  que l'on a une version locale d'applatissement pour un morphisme génériquement ouvert dans \cite{Si}. } de Hironaka (\cite{H1}) qui prouve l'existence d'une  modification $\phi:\widetilde{S}\rightarrow S$ et d'un diagramme commutatif  
$$\xymatrix{\widetilde{X}\ar[r]^{\tilde\phi}\ar[d]_{\tilde{\pi}}&{X}\ar[d]^{\pi}\\ 
 {\widetilde S}\ar[r]_{\phi}& S}$$
 dans lequel $\tilde{\pi}$ est plat, $\widetilde{X}$  constitué de toutes les composantes irréductibles du produit fibré $X\times_{S}\widetilde{S}$ qui se surjectent sur $\widetilde{S}$ (c-à-d la préimage stricte de $X$) et  $\tilde\phi$ est une modification propre. \vspace{1mm}

 \noindent Dans ce cas, étant plat, $\tilde\pi$ est géométriquement plat et on peut appliquer ce qui précède pour définir le morphisme d'image réciproque ${\mathcal L}^{k}_{S}\rightarrow \pi_{*}{\mathcal L}^{k}_{X}$.\vspace{1mm}

 \noindent 
 De plus, on vérifie aisément que ces constructions sont:
 
 \vspace{1mm}
 
 \indent $\bullet$ {\bf{invariantes par  modification propre de la base.}} En effet, soit donné un diagramme commutatif
$$\xymatrix{X_{2}\ar[r]^{\theta_{2}}\ar[d]_{\pi_{2}}&X\ar[d]^{\pi}&X_{1}\ar[l]_{\theta_{1}}\ar[d]^{\pi_{1}}\\
S_{2}\ar[r]_{\nu_{2}}&S&S_{1}\ar[l]^{\nu_{1}}}$$
dans lequel  tous les morphismes  sont propres ( ceux autres que $\pi$, $\pi_{1}$ et $\pi_{2}$  sont des modifications). Alors, les faisceaux  ${\mathcal L}^{q}_{X}$ et ${\mathcal L}^{q}_{S}$ ne dépendent pas des résolutions ou modifications choisies pour les décrire, on a
$${\mathcal L}^{k}_{X}={\theta_{2}}_{*}{\mathcal L}^{k}_{X_2}={\theta_{1}}_{*}{\mathcal L}^{k}_{X_1},\,\,{\mathcal L}^{k}_{S}={\nu_{2}}_{*}{\mathcal L}^{k}_{S_2}={\nu_{1}}_{*}{\mathcal L}^{k}_{S_1}$$
et, donc,
$$\xymatrix{{\mathcal L}^{k}_{S}\ar@{=}[r]\ar[d]&{\nu_{2}}_{*}{\mathcal L}^{k}_{S_2}\ar@{=}[r]\ar[d]&{\nu_{1}}_{*}{\mathcal L}^{k}_{S_1}\ar[d]\\
\pi_{*}{\mathcal L}^{k}_{X}\ar@{=}[r]&{\pi_2}_{*}({\theta_{2}}_{*}{\mathcal L}^{k}_{X_2})\ar@{=}[r]&{\pi_1}_{*}({\theta_{1}}_{*}{\mathcal L}^{k}_{X_1})}$$

\noindent
\vspace{1mm}

\indent $\bullet$ {\bf{compatibles avec la composition.}}\vspace{1mm}

\noindent 
Supposons donnés la composée $\xymatrix{ X\ar@/^1pc/[rr]^{h}\ar[r]_{f}& Y\ar[r]_{g}&S}$ de morphismes ouverts à fibres de dimension constante.  Comme $h$ est aussi ouvert et à fibres de dimension constante, on peut, là encore, supposer $S$ normal au vu du diagramme de changement de base donné par la normalisation 
$$\xymatrix{\tilde{X}\ar@/_2pc/[dd]_{\tilde h}\ar[r]^{\theta}\ar[d]_{\tilde{f}}&{X}\ar[d]^{f} \ar@/^2pc/[dd]^{h}\\ 
 {\tilde Y}\ar[d]_{\tilde g}\ar[r]_{\tilde\phi}& Y\ar[d]^{g}\\
 \tilde{S}\ar[r]_{\phi}&S}$$
Pour s'assurer de la relation $(\overline{f\circ g})^{*}={\bar g}^{*}\circ {\bar f}^{*}$, on peut raisonner génériquement sachant que les faisceaux sont sans torsion ou utiliser le {\bf{théorème 1}} de \cite{K3} en notant que
$$g_{*}(f_{*}{\mathcal L}^{k}_{X})\simeq g_{*}(f_{*}\theta_{*}{\mathcal L}^{k}_{\tilde X})\simeq  g_{*}({\tilde\phi}_{*}({\tilde f}_{*}{\mathcal L}^{k}_{\tilde X}))\simeq \phi_{*}({\tilde h_{*}}{\mathcal L}^{k}_{\tilde X})\simeq h_{*}{\mathcal L}^{k}_{X}$$
{\bf(ii)} Pour prouver l'existence d'une image réciproque sous cette condition d'incidence, il est possible de vérifier que les formes holomorphes génériques obtenues par image réciproque sont méromorphes et vérifient des conditions de croissance   ${\bf L}^{2}$ au sens de la \definitionref{def1}. Il est préférable, à notre sens,  d'utiliser la fonctorialité de la résolution des singularités et la caractérisation de ces formes par résolution et leur invariance par modification propre pour simplifier l'exposé. Ainsi, à tout morphisme d'espaces complexes réduits $f:X\rightarrow Y$ donné est associé un  diagramme commutatif
$$\xymatrix{\tilde{X}\ar[r]^{\tilde f}\ar[d]_{\phi'}&\tilde{Y}\ar[d]^{\phi}\\
X\ar[r]_{f}&Y}$$
dans lequel $\phi$ est une résolution des singularités de $Y$, $\tilde X$  la préimage stricte de $X$, $\phi'$ une modification propre et  $\tilde f$  la transformée stricte de $f$.\vspace{1mm}

\noindent
On pose, alors,  $${\bar f}^{*}(\xi):={\phi'}_{*}({\tilde f}^{*}(\phi^{*}\xi))$$
On peut supposer ce courant $\bar\partial$-fermé sans torsion et, donc, définissant une section du faisceau $\omega^{q}_{X}$. Mais comme $\phi'$ est une modification et que  ${\tilde f}^{*}(\phi^{*}\xi)$ est une forme holomorphe, ${\bf f}^{*}(\xi)$ 
définit bien une section du faisceau ${\mathcal L}^{q}_{X}$. \vspace{1mm}

\noindent On vérifie aisément que c'est l'unique prolongement de type ${\bf L}^{2}$ de la forme holomorphe $f^{*}(\xi\mid_{{\rm Reg}(Y)})$ sur $X\setminus f^{-1}({\rm Sing}(Y))$.\vspace{1mm}

\noindent 
La démonstration est en tout point identique à celle menée dans \cite{B5}, {\bf{(thm.4.1.1)} }rectifiée dans \cite{B6}, {\bf{(thm.1.0.1)}} et à laquelle nous renvoyons le lecteur$\blacksquare$
\end{proof}
\section{\color{blue}{Remarques.}}
On peut remarquer que ces constructions supposent la surjectivité des morphismes ou la condition d'incidence qui permettent des raisonnements sur les parties génériques. Comme on l'a vu dans \cite{K3}, ces conditions ne sont pas suffisantes pour déterminer une image réciproque pour les faisceaux $\omega^{bullet}_{Z}$. Il suffit de le constater sur l'exemple simple de la normalisation
$\nu:{\Bbb C}^{2}\rightarrow Z:=\{(x,y,z)\in {\Bbb C}^{3}: xy^{2}=z^{2}\}$ envoyant $(u,v)$ sur $(u^{2}, v, uv)$ et la forme méromorphe $\displaystyle{\frac{dz}{y}}$.

\phantomsection\addcontentsline{toc}{part}{Le \theoremref{T1} pour un morphisme arbitraire.}
%
%
\svnid{$Id: S07-proof-setup.tex 269 2020-01-20 11:28:53Z kebekus $}

\section{\color{blue}{Morphisme de restriction: cas particulier.}}
\subversionInfo

\approvals{Mohamed & yes}
\par\vspace{2mm}

\subsection{{La  restriction sur les faisceaux ${\mathcal L}^{\bullet}_{X}$}.}
Ayant pour objectif de définir l'image réciproque par un morphisme arbitraire sur ces faisceaux et le problème étant de nature local, on se ramène fondamentalement au cas où ce morphisme est la composé d'un plongement suivi d'une projection (factorisation par le graphe) et, donc, au cas d'un plongement, celui de la projection étant trivial!
\vspace{1mm}

\noindent
La partie substancielle du \theoremref{T1} consiste à définir cette restriction car le reste n'en est qu'une conséquence. Dans ce paragraphe, nous montrons le   
\Th{}{}\label{T4}   Soit $X$ un espace analytique complexe réduit de dimension pure $m$. Alors, pour tout sous espace fermé $Y$ de $X$ muni d'un plongement local $i$, il existe un morphisme ${\mathcal O}_{Y}$-linéaire  de faisceaux cohérents $ {\mathfrak i}^{*}:{\mathcal 
L}^{q}_{X}\rightarrow i_{*}{\mathcal L}^{q}_{X}$
prolongeant  naturellement  la restriction usuelle des formes holomorphes rendant (modulo torsion du faisceau des formes holomorphes)   commutatif le diagramme 
$$\xymatrix{\Omega^{q}_{X}\ar[r]\ar[d]& i_{*}(\Omega^{q}_{Y})\ar[d]\\
{\mathcal L}^{q}_{X}\ar[r]_{{\mathfrak i}^{*}}& i_{*}{\mathcal L}^{q}_{Y}}$$
 \noindent
De plus, on a, pour toute forme holomorphe $\alpha$ et toute section $\xi$ du faisceau 
${\mathcal L}^{q}_{X}$, la relation de compatibilité avec le produit extérieur
$$ {\mathfrak i}^{*}(\alpha \wedge \xi) =
i^{*}(\alpha)\wedge {\mathfrak i}^{*}(\xi)$$
et   avec la composition des plongements locaux $\xymatrix{ Y'\ar@/^1pc/[rr]^{i"}\ar[r]_{i'}& Y\ar[r]_{i}&X}$ 
$${\mathfrak i"}^{*}=  {\mathfrak i'}^{*}\circ{\mathfrak i}^{*}$$\rm
 \vspace{1mm}
\rm
 \noindent 
\subsection{Construction et définition d'une restriction dans le cas lisse.}\vspace{1mm}

\noindent
 La définition que nous avons en vue utilise l'intégration sur les fibres. Avant de la proposer dans le cas général, nous l'éprouvons dans le cas lisse et montrons qu'elle s'identifie à la restriction  usuelle.  Nous nous attarderons pas sur le caractère canonique de la construction qui sera étudié en détail dans la preuve générale.\vspace{1mm}
 
 \noindent 
Rappelons qu'une  forme de K\"ahler relative  pour un morphisme propre et surjectif (donc localement k\"ahlérien) ) $\pi:X\rightarrow T$ représente une
 classe de  $H^{1}(X, \Omega^{1}_{X})$ qui est en fait 
une section de 
$H^{0}(T, {\rm I}\!{\rm R}^{1}\pi_{*}\Omega^{1}_{X/T})$ dont la restriction
à chaque fibre est une classe de K\"ahler;  nous renvoyons le lecteur à \cite{Fu-Sc}, \cite{Sc} et \cite{Va} pour la  définition des formes de K\"ahler 
relatives.\vspace{1mm}

\noindent

\Prop{}{}\label{lem1} Soit $\pi: X\rightarrow Y$ un morphisme propre et surjectif d'espaces complexes réduits dont la dimension générique des fibres est $n$. Soient $q$ un entier naturel,  $\xi$ une section du faisceau ${\mathcal L}^{q}_{X}$ et $\omega$ une forme de K\"ahler relative. Alors, la forme $\int_{X/Y}\xi\wedge\omega^{n}$, définie génériquement sur $Y$, se prolonge en une unique section du faisceau ${\mathcal L}^{q}_{Y}$.\rm
\begin{proof} Comme $\pi$ est propre et surjectif,  l'image directe $\pi_{*}(\xi\wedge \omega^{n})$ définit un courant de type $(q,0)$ et, donc, modulo torsion, une section du faisceau $\omega^{q}_{Y}$.  Mais $Y$ étant réduit, il existe un ouvert dense de $Y$ sur lequel  $\pi$ est géométriquement plat (resp. plat).  Alors, le théorème d'applatissement géométrique de  (\cite{B2}, \cite{Si}) (resp. d'applatissement algébrique de  \cite{H1})   montre qu'il existe  une modification $\nu:\tilde{Y}\rightarrow Y$ et un diagramme commutatif
 $$\xymatrix{\tilde{X}\ar[r]^{\Theta}\ar[d]_{\tilde\pi}&X\ar[d]^{\pi}\\
\tilde{Y}\ar[r]_{\nu}&Y}$$ 
dans lequel $\tilde\pi$ est géométriquement plat (resp. plat) et $\Theta$ une modification propre.  Le théorème d'intégration sur les fibres d'un morphisme géométriquement plat (cf \cite{K0}) assure que l'intégrale sur les fibres $\int_{\tilde{X}/\tilde{Y}}\Theta^{*}(\xi\wedge\omega^{n})$  définit une section du faisceau ${\mathcal L}^{q}_{\tilde{Y}}$ (signalons au passage que $\Theta^{*}\omega$ est une forme de Kähler relative) que l'on peut aussi voir  comme  courant image directe  ${\tilde\pi}_{*}(\Theta^{*}(\xi\wedge \omega^{n}))$.   Or,   les  courants $\bar\partial$-fermés   $\nu_{*}(\int_{\tilde{X}/\tilde{Y}}\Theta^{*}(\xi\wedge\omega^{n}))$ et  $\pi_{*}(\xi\wedge \omega^{n})$ coincident   génériquement  sur $Y$ (par construction et en vertu de la formule du changement de  base pour un morphisme propre!). Par conséquent, modulo torsion et, donc, en tant que sections du faisceau $\omega^{q}_{Y}$, ils coincident partout sur $Y$. Comme  on a   ${\mathcal L}^{q}_{Y}=\nu_{*}{\mathcal L}^{q}_{\tilde{Y}}$, on en déduit le résultat. \vspace{1mm}

\noindent En termes de courants, on peut remarquer que l'égalité
$$\langle{\tilde\pi}_{*}(\Theta^{*}(\xi\wedge \omega^{n})), \phi\rangle= \int_{\tilde X}\Theta^{*}(\xi\wedge \omega^{n})\wedge {\tilde\pi}^{*}(\phi) $$
donne,en particulier,  pour  $\phi=\nu^{*}(\psi)$,  
$$\langle{\tilde\pi}_{*}(\Theta^{*}(\xi\wedge \omega^{n})), \phi\rangle= \langle{\pi}_{*}(\xi\wedge \omega^{n}), \psi\rangle$$
qui suffit aussi  pour conclure$\,\blacksquare$
\end{proof}\vspace{1mm}

\Prop{}{}\label{P'1} Soient $Z$ une variété analytique complexe de dimension pure $m$ et  $T$ un sous ensemble analytique de codimension $r$ dans  $Z$ muni d'un plongement local $i$. Soient  $\pi:{\tilde Z}\rightarrow Z$  une modification propre de centre lisse  $T$ et  $\omega$ une forme  de K\"ahler $\pi$-relative.  Soit $Y$ un sous espace de $Z$ muni d'un plongement local $i$.  Alors,  pour toute section $\xi$ du faisceau $\Omega^{q}_{Z}$ sur un ouvert $U$, on pose $${\mathcal R}^{\pi,\omega}_{Z,Y}(\xi)={\tilde \pi}_{*}(\pi^{*}(\xi)\wedge \omega^{k}|_{\tilde Y})\,\,\rm {avec}$$
 \indent $\bullet$ $k=0$  si ${Y\cap T}=\emptyset$ (resp. ${Y\cap T}\subsetneq Y$ )  et  $\tilde\pi$ est l'isomorphisme (resp. la modification propre) induit (e) par $\pi$ sur  la préimage totale (resp.  stricte) ${\tilde Y}:=\pi^{-1}(Y)$ sur $Y$,\vspace{1mm}

\indent
$\bullet$ $k=r-1$  si $Y\subseteq T$ et  ${\tilde \pi}$ le morphisme induit par $\pi$ sur la préimage totale de  $Y$.\vspace{1mm}

\noindent
Alors,   ${\mathcal R}^{\pi,\omega}_{Z,Y}(\xi)$ est un courant $\bar\partial$-fermé sans torsion définit une section du faisceau ${\Omega}^{q}_{Y}$ indépendante des choix faits et induit, de ce fait, un morphisme ${\mathcal O}_{Y}$-linéaire  de faisceaux cohérents $ {\mathfrak i}^{*}:{\Omega}^{q}_{Z}\rightarrow i_{*}{\Omega}^{q}_{Y}$ coincidant avec le morphisme de restriction usuel.\rm\vspace{1mm}

\noindent
\begin{proof} Soient $Z$ et $T$ comme dans l'énoncé et  $Y$ un sous espace analytique de $Z$ comme l'on peut supposer être un sous ensemble analytique puisque cette construction est nature locale.  Soient $\pi:{\tilde Z}\rightarrow Z$ une modification propre de centre $T$ et $\omega$ une forme de K\"ahler relative pour $\pi$. Alors, pour tout ouvert $U$ de $Z$ et toute section $\xi$ du faisceau $  \Omega^{q}_{Z}$ sur $U$, on pose $${\mathcal R}^{\pi,\omega}_{Z,Y}(\xi)={\tilde \pi}_{*}(\pi^{*}(\xi)\wedge \omega^{k}|_{\tilde Y})$$
avec: \vspace{1mm}

\indent $\bullet$ $k=0$  si ${Y\cap T}=\emptyset$ (resp. ${Y\cap T}\subsetneq Y$ )  et  $\tilde\pi$ est l'isomorphisme (resp. la modification propre) induit (e) par $\pi$ sur  la préimage totale (resp.  stricte) ${\tilde Y}:=\pi^{-1}(Y)$ sur $Y$,\vspace{1mm}

\indent
$\bullet$ $k=r-1$  si $Y\subseteq T$ et  ${\tilde \pi}$ le morphisme induit par $\pi$ sur la préimage totale de  $Y$.\vspace{2mm}

\noindent
{\bf(a)} $\star$ Si ${Y\cap T}=\emptyset$, le résultat est trivial puisque $\pi$ est un isomorphisme analytique de ${\tilde Z}\setminus \pi^{-1}(T)$ sur $Z\setminus T$.  Dans ce cas, il est évident que le courant ${\tilde\pi}_{*}((\pi^{*}(\xi))|_{\tilde Y})$ coincide naturellement avec la  restriction usuelle de la forme $\xi$ puisque la préimage totale  ${\tilde Y}:= \pi^{-1}(Y)$ est isomorphe à $Y$. \par
$\star$ Si l' intersection ${Y\cap T}$ est d'intérieur vide dans $Y$,  le morphisme $\pi$ induit sur la préimage stricte ${\tilde Y}$ (qui est l'adhérence de  $\pi^{-1}(Y\setminus {Y\cap T})$ dans le produit fibré ${\tilde Z}\times_{Z}Y$) une modification propre de centre ${Y\cap T}$ installée  dans le diagramme commutatif  
$$\xymatrix{\tilde{Y}\ar[d]_{\tilde \pi}\ar[r]^{\tilde i}&{\tilde Z}\ar[d]^{\pi}\\
Y\ar[r]_{i}&Z}$$
Alors, pour toute section $\xi$ du faisceau ${\Omega}^{q}_{Z}$, le courant ${\mathcal R}^{\pi,\omega}_{Z,Y}(\xi):={\tilde\pi}_{*}({\tilde i}^{*}(\pi^{*}(\xi)))$ est $\bar\partial$-fermé et sans torsion puisque non trivial sur $Y\setminus Y\cap T$ car ce n'est rien d'autre que la restriction usuelle de la forme holomorphe $\xi$ d'après ce que l'on a dit précédemment. Il définit, donc, naturellement une section du faisceau $\omega^q_{Y}$ qui, d'après le  lemme de Dolbeault-Grothendieck, est une section du faisceau $\Omega^q_{Y}$ sur la partie lisse de $Y$. En fait,  la formule de projection et la  définition des courants d'intégration sur $\tilde{Y}$ et $Y$ respectivement, nous donnent
$$i_{*}{\tilde \pi}_{*}(\pi^{*}(\xi)|_{\tilde Y})= {\pi}_{*}({\tilde i}_{*}({\tilde i}^{*}(\pi^{*}(\xi))= {\pi}_{*}(\pi^{*}(\xi)\wedge [\tilde Y])=\xi\wedge \pi_{*}([\tilde Y])=i_{*}i^{*}(\xi)$$
qui, au vu de l'égalité $\pi_{*}([\tilde Y])= [ Y]$ et de  l'injectivité du  foncteur $i_{*}$, nous donne
$${\mathcal R}^{\pi,\omega}_{Z,Y}(\xi)=i^{*}(\xi)$$
{\bf (b)} Si $Y\subseteq T$,  on a  le diagramme commutatif
$$\xymatrix{\pi^{-1}(Y)\ar[d]_{\tilde \pi}\ar[r]^{\tilde i}& \pi^{-1}(T)\ar[r]\ar[d]^{\pi'}& {\tilde Z}\ar[d]^{\pi}\\
Y\ar[r]_{i}&T\ar[r]&Z}$$
dans lequel les morphismes  $\pi'$ et $\tilde\pi$ sont lisses et $(r-1)$-équidimensionnels; le premier par construction puisque  $\pi^{-1}(T)$ s'identifie au projectifié du fibré normal de $T$ dans $Z$, ${\Bbb P}({\mathcal N}_{T/Z})$,  le second par stabilité par changement de base  des notions de  lissité et de  la constance de la dimension des fibres. Il apparait clairement qu'il nous suffit de traiter le cas où $Y=T$ (le cas général s'y ramène  puisque ce courant image directe coincide avec l'intégration sur les fibres et que cette opération est compatible aux changements de bases. Alors, choisissant $\omega$ telle $\int_{{\tilde Y}/Y}\omega^{r-1}=1$,  on pose $${\mathcal R}^{\pi,\omega}_{Z,Y}(\xi):={\tilde \pi}_{*}(\pi^{*}(\xi)\wedge \omega^{r-1}|_{\tilde Y})$$
Or, on a, au sens des courants,
$${\mathcal R}^{\pi,\omega}_{Z,Y}(\xi):={\tilde \pi}_{*}(\pi^{*}(\xi)\wedge \omega^{r-1}|_{\tilde Y})={\tilde \pi}_{*}({\tilde \pi}^{*}(i^{*}\xi)\wedge \omega^{r-1}|_{\tilde Y})$$
qui nous donne, grâce à la formule de projection et la normalisation choisie
$${\mathcal R}^{\pi,\omega}_{Z,Y}(\xi)=i^{*}(\xi)$$
Remarquons, au passage, que l'on a de façon équivalente (en utilisant les couranst d'intégration et l'injection de l'image directe $i_*$)
$$i_{*}{\tilde \pi}_{*}(\pi^{*}(\xi)\wedge \omega^{r-1}|_{\tilde Y})= {\pi}_{*}({\tilde i}_{*}({\tilde i}^{*}(\pi^{*}(\xi)\wedge \omega^{r-1}))= {\pi}_{*}(\pi^{*}(\xi)\wedge \omega^{r-1}\wedge [\tilde Y])=\xi\wedge \pi_{*}(\omega^{r-1}\wedge [\tilde Y])=\xi$$
L'égalité  $\pi_{*}(\omega^{r-1}\wedge [\tilde Y])= [ Y]$ s'obtient grâce à la normalisation choisie et au fait que  ceux sont deux courants $\bar\partial$-fermés sans torsion, donc définissant des sections du faisceau $\omega^{q}_{Y}$, coincidant génériquement sur $Y$ (par exemple, aux points en lesquels $\tilde\pi$ est submersif). \vspace{1mm}

\noindent
On admettra pour l'instant, les propriétés d'invariance de la construction et son caractère intrinsèque qui seront abordés dans toute leur généralité et de façon détaillée dans la preuve du cas général.$\,\blacksquare$

\end{proof} \vspace{2mm}

\noindent
\section{\color{blue}{Morphisme de restriction: cas général}}

\Prop{}{}\label{P'2} Soit $\pi: Z'\rightarrow Z$ un morphisme localement k\"ahlérien propre et  surjectif d'espaces complexes réduits avec $Z'$ lisse et doté d'une forme de Kähler relative $\omega$. Soit $Y$ (resp. $Y'$)   un  sous espace réduit  de codimension $r$ dans  $Z$ (resp. $Z'$)  muni d'un morphisme $\pi': Y'\rightarrow Y$ propre surjectif génériquement $n$-équidensionnel tel que  le diagramme  $$\xymatrix{{Y'}\ar[d]_{\pi'}\ar[r]^{i'}&{Z'}\ar[d]^{\pi}\\
 Y\ar[r]^{i}& Z}$$
soit commutatif. On  choisira  $\omega$ de sorte à ce que $\int_{Y'/Y}\omega^{n}=1$.\vspace{1mm}

\noindent Alors,   pour toute section $\xi$ du faisceau  ${\mathcal L}^{q}_{Z}$,  le courant défini  par   ${\mathcal R}^{\pi,\omega}_{Z,Y}(\xi):={\pi'_{*}}(\pi^{*}(\xi)\wedge \omega^{n}|_{\tilde Y})$ est $\bar\partial$-fermé et définit une unique section de  ${\mathcal L}^{q}_{Y}$.\vspace{1mm}

\noindent  
De plus, cette correspondance, induit un morphisme ${\mathcal O}_{Y}$-linéaire ${\mathcal L}^{q}_{Z}\rightarrow i_{*}{\mathcal L}^{q}_{Y}$ prolongeant  la restriction usuelle des formes holomorphes et rendant commutatif le diagramme 
$$\xymatrix{\Omega^{q}_{Z}/{\mathcal T}_{Z}\ar[r]\ar@{^{(}->}[d]& i_{*}(\Omega^{q}_{Y}/{\mathcal T}_{Y})\ar@{^{(}->}[d]\\
{\mathcal L}^{q}_{Z}\ar[r]_{{\mathcal R}^{\pi,\omega}_{Z,Y}}& i_{*}{\mathcal L}^{q}_{Y}}$$
où ${\mathcal T}_{Z}$ (resp. ${\mathcal T}_{Y}$) désigne le sous faisceau de torsion de $\Omega^{q}_{Z}$ (resp. $\Omega^{q}_{Y}$). \rm
\begin{proof} \par\noindent
{\bf(a)} Tout comme dans la \propositionref{P'1}, on a, pour toute forme holomorphe $\xi$ sur $Z$,  
$${\mathcal R}^{\pi,\omega}_{Z,Y}(\xi)=i^{*}(\xi)$$
 qui découle simplement de la formule de projection et de la normalisation choisie $\pi'_{*} \omega^{n}|_{Y'}=1$ puisque 
$${\mathcal R}^{\pi,\omega}_{Z,Y}(\xi):=\pi'_{*}(\pi^{*}(\xi)|_{Y'}\wedge \omega^{n}|_{Y'})=\pi'_{*}(\pi'^{*}(\xi|_{Y})\wedge \omega^{n}|_{Y'})=(\xi|_{Y})\pi'_{*} \omega^{n}|_{Y'} $$
On définit ainsi un morphisme ${\mathcal R}^{\pi,\omega}_{Z,Y}: \Omega^{q}_{Z}\rightarrow i_{*}\Omega^{q}_{Y}$ dont on déduit naturellement (par passage au quotient)
$ {\mathcal R}^{\pi,\omega}_{Z,Y}:\Omega^{q}_{Z}/{\mathcal T}_{Z}\rightarrow i_{*}(\Omega^{q}_{Y}/{\mathcal T}_{Y}) $  puisque, de façon générale, un morphisme de faisceaux analytiques cohérents sur un espace complexe réduit applique toujours le sous faisceau de torsion de l'un sur le sous faisceau de torsion de l'autre. Alors, en supposant vrai le fait que ${\mathcal R}^{\pi,\omega}_{Z,Y}(\xi)$ définisse une section du faisceau  ${\mathcal L}^{q}_{Y}$, l'existence et la  commutativité du diagramme reliant la restriction usuelle et ${\mathcal R}^{\pi,\omega}_{Z,Y}$ est évidente puisque l'injections naturelle
$\Omega^{q}_{\tilde Y}/{{\mathcal T}_{\tilde Y}}\rightarrow {\mathcal L}^{q}_{\tilde Y}$ nous donne  le  morphisme injectif ${\tilde\pi}_{*}(\Omega^{q}_{\tilde Y}/{{\mathcal T}_{Y}})\rightarrow {\mathcal L}^{q}_{Y}$ (puisque ${\tilde\pi}_{*}{\mathcal L}^{q}_{\tilde Y}={\mathcal L}^{q}_{Y}$) et de même sur $Z$.
 \vspace{1mm}

\noindent
{\bf(b)} Il nous faut montrer que, pour toute section  $\xi$ du faisceau ${\mathcal L}^{q}_{Z}$, on définit bien par ce procédé un section du faisceau ${\mathcal L}^{q}_{Y}$. \vspace{1mm}

\noindent
Comme le courant ${\mathcal R}^{\pi,\omega}_{Z,Y}(\xi):=\pi'_{*}(\pi^{*}(\xi)|_{Y'}\wedge \omega^{n}|_{Y'})$ est  de type $(q,0)$ et $\bar\partial$-fermé sur $Y$ définissant, sur sa partie lisse, en vertu du lemme de Dolbeault-Grothendieck,  une $q$-forme holomorphe d'ailleurs non triviale (comme on peut s'en convaincre aisément en nous restreignant sur la partie régulière du morphisme $\pi'$), il est sans torsion et définit, par suite, une unique section du faisceau $\omega^{q}_{Y}$.  Pour montrer que cette construction nous donne, en fait, une section du faisceau ${\mathcal L}^{q}_{Y}$, on peut  soit  utiliser le fait que, pour toute désingularisation de $Y$, ce procédé nous fournit une forme holomorphe et, par conséquent, ${\mathcal R}^{\pi,\omega}_{Z,Y}(\xi)$ est en fait une section du faisceau  ${\mathcal L}^{q}_{Y}$ ( cet argument fait appel à une formule de transitivité que l'on verra un peu plus loin), soit  la \propositionref{lem1} qui montre que le courant obtenu est bien une section de ce faisceau. En effet, comme  cette image directe  est donnée, au moins génériquement, par l'intégrale sur les fibres  $\int_{Y'/Y}\pi^{*}(\xi)|_{Y'}\wedge \omega^{n}|_{Y'}$, on montre, alors, que l'on construit, ainsi, une forme holomorphe sur $Y$ qui n'est autre    que la restriction usuelle de la forme $\xi$ ce dont on s'en convainc aisément, vu la nature locale de notre problème est de nature locale sur $Z$ (et $Y$). On peut, quitte à rétrécir les données, le supposer plonger dans un ouvert de Stein $W$ d'un certain espace numérique ${\Bbb C}^{N}$ de sorte à ce que  toute forme holomorphe sur $Z$ soit induite (par passage au quotient) par une forme holomorphe sur $W$ et, par suite, supposer $Z$ lisse. En effet, dans ce cadre local,  il existe une désingularisation plongée de $Z$ donnée par  l'installation 
$$\xymatrix{{Y'}\ar[d]_{\pi'}\ar[r]^{i'}&\tilde{Z}\ar[d]_{\pi}\ar[r]^{\sigma'}&W'\ar[d]^{\phi}\\
 Y\ar[r]^{i}&Z\ar[r]^{\sigma}& W}$$
 On se ramène ainsi au cas lisse traité dans la \propositionref{P'1} puisque, pour toute forme holomorphe $\xi$ sur $W$, on a, en posant $\tau:=i\circ \sigma$ (resp. $\tau':=i'\circ \sigma'$  la relation 
$$\tau^{*}(\xi)=\pi'_{*}(\phi^{*}(\xi)|_{Y'}\wedge \omega^{n}|_{Y'})$$
En effet, on a, au sens des courants,
$$\tau_{*}\pi'_{*}(\phi^{*}(\xi)|_{Y'}\wedge \omega^{n}|_{Y'})=\phi_{*}(\tau'_{*}\tau'^{*}(\phi^{*}(\xi)\wedge\omega^{n}))=\phi_{*}(\phi^{*}(\xi)\wedge\omega^{n}\wedge[Y'])=\xi\wedge \phi_{*}(\omega^{n}\wedge[Y'])$$
Mais,  $\phi_{*}(\omega\wedge[Y'])=[Y]$ comme on s'en convainc aisément grâce à la submersivité de $\pi'$. Ainsi, obtient-on l'égalité des courants
$$\tau_{*}\tau^{*}\xi=\tau_{*}\pi'_{*}(\phi^{*}(\xi)|_{Y'}\wedge \omega^{n}|_{Y'})$$
garantissant l'égalité cherchée par injectivité de l'image directe.\vspace{1mm}

\noindent
Remarquons que  si le sous ensemble au dessus duquel $\pi'$ est submersif est de codimension deux (ou plus), toute forme holomorphe de $Y$ définie en dehors de sous ensemble est automatiquement une section du faisceau  ${\omega}^{q}_{Y}$ qui est de profondeur au moins deux.\vspace{1mm}

\noindent Cette opération  d'intégration le long des fibres applique le faisceau de torsion ${\mathcal T}_{Z}$ sur le faisceau de torsion ${\mathcal T}_{Y}$; d'ailleurs, le lecteur désireux de voir les détails de cette assertion peut se référer  à (\cite{B3}, {\it prop (1.0.1)}, {\it cor  (1.0.2)}; p.3-4)). La \propositionref{P'2} permet de définir un morphisme de faisceaux ${\mathcal R}^{\pi,\omega}_{Z,Y}$ entre les faisceaux  ${\mathcal L}^{q}_{Z}$ et $i_{*}{\mathcal L}^{q}_{Y}$. Comme le faisceaux cohérent $\Omega^{q}_{Z}/{\mathcal T}_{Z}$ (resp. $ i_{*}(\Omega^{q}_{Y}/{\mathcal T}_{Y}$) s'injecte dans ${\mathcal L}^{q}_{Z}$ (resp. $i_{*}{\mathcal L}^{q}_{Y}$), la commutativité du diagramme de la proposition  devient évidente.  En effet, partant d'une forme holomorphe  sans torsion $\xi$,  on obtient une forme holomorphe $i^{*}(\xi)$  sans torsion sur $Y$ pouvant être vu comme une section du faisceau ${\mathcal L}^{q}_{Y}$. Mais alors, en tant que section de ce dernier, elle coincide génériquement sur $Y$ ( en les points de submersivité de $\pi'$ par exemple) avec $\int_{Y'/Y}\pi^{*}(\xi)|_{Y'}\wedge\omega^{n}|_{Y'}$,  d'après {\bf(a)}.  Comme le faisceau ${\mathcal L}^{q}_{Y}$  est sans torsion, ces deux sections s'identifient sur $Y$ tout entier$\,\blacksquare$

\end{proof} 
\section{\color{blue}{ La preuve du \theoremref{T4}.}}
\par

\subsubsection{\bf{Construction et définition de la restriction}.}\par\noindent
La définition que l'on peut donner nécessite de prendre en considération l'incidence de $Y$ avec le lieu singulier  de $Z$.  Le cas le plus emblématique est bien entendu celui où $Y$ est contenu dans le lieu singulier; le cas où $Y$ est en dehors du lieu singulier de $Z$ ne présente aucune difficulté et a déjà été  traité, en grande partie,  dans la \propositionref{P'1} et la \propositionref{P'2} . \vspace{1mm}

\noindent
Soit $Z$ un espace analytique réduit et $Y$ un sous espace réduit de $Z$ et   
$$\xymatrix{Y'\ar[r]^{i'}\ar[d]_{\pi'}&Z'\ar[d]^{\pi}\\
Y\ar[r]_{i}&Z}$$ 
un diagramme commutatif dans lequel $Z'$ est une variété lisse, $\pi$  un morphisme localement k\"ahlérien, $i$ et $i'$ des plongements. Soit $\omega$ une forme de K\"ahler relative. Alors, on pose, pour toute section $\xi$ du faisceau ${\mathcal L}^{q}_{Z}$, 
\[{\mathcal R}^{\pi,\omega}_{Z,Y}(\xi):={\pi'_{*}}(\pi^{*}(\xi)|_{Y'} \wedge \omega^{k}|_{Y'})\]
avec $k=0$ si $\pi'$ est une modification propre et $k=n$ si $\pi'$ est propre génériquement  $n$-équidimensionnel.\vspace{1mm}

\noindent 
En anticipant un peu sur l'exposé, on ne va considérer que les morphismes de désingularisations car, comme nous le verrons,  la construction est indépendante du morphisme localement k\"ahlérien choisi.  Ces deux expressions correspondent aux cas de figure suivants:\vspace{1mm}

\indent 
{\bf (a)  $Y$ non entièrement contenu dans le lieu singulier de $Z$.}\par\noindent
Considérons le diagramme commutatif suivant
$$\xymatrix{{Y'}\ar[rd]_{\pi'}\ar@/^2pc/[rr]^{i'}\ar@{^{(}->}[r]^{i'}&{Z'}\times_{Z}Y\ar[d]^{p_{2}}\ar[r]^{p_{1}}&{Z'}\ar[d]^{\pi}\\
&Y\ar[r]_{i}&Z}$$
dans lequel $\pi$ est un morphisme de désingularisation,  $p_{1}$ (resp. $p_{2}$) sont les projections canoniques,  $Y'$ est la préimage stricte décrite par l'adhérence de ${p_{2}}^{-1}(Y\setminus{\rm Sing}(Z)\cap Y)$ dans le produit fibré.\par\noindent
Il est bien connu que le  morphisme $\pi'$, induit par la projection $p_{2}$, est  une modification propre. Dans ce cas, si $\xi$ est une section du faisceau ${\mathcal L}^{q}_{Z}$, 
$${\mathcal R}^{\pi}_{Z,Y}(\xi):={\pi'}_{*}({ i'}^{*}(\pi^{*}(\xi)))$$
Comme $Y$ n'est pas entièrement contenu dans ${\rm Sing}(Z)$, $Y_{0}:=Y\cap {\rm Reg}(Z)$ est un ouvert dense de $Y$, la forme holomorphe ${ i'}^{*}(\pi^{*}(\xi))$ n'est ni triviale ni de torsion  sur ${ Y'}$. Dans ce cas, 
 le courant de type $(q,0)$ et $\bar\partial$-fermé ${\mathcal R}^{\pi}_{Z,Y}(\xi)$ définit une section du faisceau ${\mathcal L}^{q}_{Y}:={\pi'_{*}}{\mathcal L}^{q}_{Y'}$ non triviale puisque l'on sait, d'après  la \propositionref{P'1} appliquée à ${\rm Reg}(Z)$ et $Y_0$, que  la $q$-forme holomorphe ${\mathcal R}^{\pi}_{{\rm Reg}(Z),Y_{0}}(\xi\mid_{{\rm Reg}(Z)})$, qui n'est rien d'autre que la restriction de la forme holomorphe $\xi\mid_{{\rm Reg}(Z)}$ à $Y_0$, est non triviale. Enfin, la \propositionref{P'2} montre qu'au niveau des formes holomorphes cette opération coincide avec la restriction usuelle c'est-à-dire que ${\mathcal R}^{\pi}_{Z,Y}(\xi)$ s'identifie à la restriction de la forme holomorphe  $\xi$ à $Y$.  Ces remarques suffisent à construire le diagramme commutatif naturel reliant la restriction usuelle sur les formes holomorphes et ces formes méromorphes.
 \vspace{1mm}

\indent 
{\bf (b)  $Y$ est contenu dans le lieu singulier de $Z$.}\vspace{1mm}

\noindent  La stratégie consiste à trouver un ouvert dense $Y_{0}$ dans $Y$ sur lequel $\pi'$ est géométriquement plat pour pouvoir appliquer le  théorème d'intégration sur les fibres de \cite{K0}  garantissant  l'existence d'un morphisme 
${\rm I}\!{\rm R}^{n}\pi_{*}\Omega^{n+q}_{Y'_{0}}\rightarrow{\mathcal L}^{q}_{Y_{0}}$. On constate ensuite que les 
  sections obtenues génériquement sur $Y$ se 
prolongent en sections globales du faisceau ${\mathcal L}^{q}_{Y}$.\vspace{1mm}

\noindent 
Soit $\pi: Z'\rightarrow Z$ un morphisme localement k\"ahlérien et surjectif d'une 
variété analytique complexe $Z'$ sur $Z$. Soit $Y$ un sous espace réduit de $Z$ et $Y':= \pi^{-1}(Y)$ sa préimage totale.  On note  $\pi'$ la 
restriction de $\pi$ 
à $Y'$ et $n := 
{\rm dim}Y' - {\rm dim} Y$, la dimension générique des fibres de ce 
morphisme. Soit $\omega$ une forme de K\"ahler  relative au morphisme $\pi'$ telle que  $\displaystyle{\int_{Y'/Y}\omega^{n}:=\int_{\pi^{-1}(y)}\omega^{n} = 1}$, $\forall y\in Y$.
\par\noindent Alors, pour tout entier $q$  et tout ouvert $U$ de Stein 
dans $Z$,  on définit  un morphisme 
  $$ {\mathcal R}^{\pi,\omega}_{Z,Y}: \Gamma(U, {\mathcal 
L}^{q}_{Z})\rightarrow 
\Gamma( U\cap Y,{\mathcal L}^{q}_{Y})$$ 
en posant
$${\mathcal R}^{\pi,\omega}_{Z,Y}(\xi):= 
\int_{Y'/Y}\pi^{*}(\xi)|_{Y'} \wedge \omega^{n}|_{Y'}$$  
\noindent
dépendant, à priori, du morphisme $\pi$ et de la forme de K\"ahler $\omega$.\par\noindent
On peut se convaincre facilement que l'on peut supposer $\pi'$ géométriquement plat sans enfreindre la généralité. En effet, il est clair que l'on peut supposer $Y'$ irréductible (on peut d'ailleurs utiliser un théorème du type  Bertini (cf \cite{Ve}, \cite{F2}) pour nous ramener au
  cas où $Y'$ est lisse). Le morphisme $\pi'$ étant propre et $Y'$ irréductible, le théorème de Cartan-Remmert sur la semi continuité de la dimension des fibres d'un morphisme propre ( \cite{Fi}, p.140) garantit la son $n$-équidimensionnalité en dehors d 'un fermé $Y_{0}$ de  codimension au moins deux dans $Y$. Comme le courant ${\mathcal R}^{\pi,\omega}_{Z,Y}(\xi)$ de type $(q,0)$ $\bar\partial$-fermé est sans torsion puisque non trivial sur la partie régulière du morphisme $\pi'$ (cf \propositionref{P'2}), il définit une unique section du faisceau $\omega^{q}_{Y}$. Or, ce dernier étant de profondeur au moins deux,  on peut supposer $\pi'$ $n$-équidimensionnel. Enfin, comme les faisceaux ${\mathcal L}^{q}_{Y}$ sont invariants par modification propre, on peut supposer $Y$ normal (après avoir effectué le changement de base par la normalisation). Mais, un morphisme équidimensionnel sur une base normale est géométriquement plat. D'où la réduction annoncée. \vspace{1mm}

\noindent Sous cette hypothèse, on déduit  du \corollaryref{C2}, que le morphisme 
$${\rm I}\!{\rm R}^{n}\pi'_{*}\Omega^{n+q}_{Y'}\rightarrow {\mathcal L}^{q}_{Y}$$
est naturellement déduit de l'intégration de classes de cohomologie (cf \cite{K0})
$${\rm H}^{n}(Y', \Omega^{n+q}_{Y'})\rightarrow \Gamma(Y, {\mathcal L}^{q}_{Y})$$
Comme $\pi^{*}(\xi)|_{Y'} \wedge \omega^{n}|_{Y'}$ représente naturellement une classe de cohomologie ${\rm H}^{n}(Y', \Omega^{n+q}_{Y'})$, on voit que cette image directe de courants s'interprète en termes d'intégration de classes de cohomologie le long des fibres du morphisme géométriquement plat $\pi'$.\vspace{1mm}

\noindent 
Comme nous le montre la \propositionref{P'2},  ce procédé fournit bien des section du faisceau  ${\mathcal L}^{q}_{Y}$.
 \vspace{1mm}

 \noindent
{\bf 3.4.2. Propriétés d'invariance de cette définition.}\par
 {\bf(a) l'indépendance de la construction vis-à-vis de la forme de K\"ahler choisie.}\par\noindent
Soit
 
$$\xymatrix{Y'\ar[d]_{\pi'}\ar[r]^{i'}&Z'\ar[d]^{\pi}\\
 Y\ar[r]_{i}&Z}$$
\noindent le diagramme commutatif ( qui est un carré cartésien) dans 
lequel $i$ et $i'$ sont les inclusions naturelles.\par\noindent
Considérons deux formes de K\"ahler relatives $\omega_{1}$ et $\omega_{2}$    vérifiant: 
$\displaystyle{\int_{Y'/Y}\omega_{1}^{n} = \int_{Y'/Y}\omega_{2}^{n} = 
1}$ (rappelons que cela s'écrit aussi $\displaystyle{\int_{\pi^{-1}(y)}\omega_{1}^{n} = \int_{\pi^{-1}(y)}\omega_{2}^{n} = 
1}$ $\forall y\in Y$), alors $$[\omega_{1}]^{n}=[\omega_{2}]^{n}$$
En effet,  le théorème de Lelong (\cite{Le}) et le lemme du découpage permettent de supposer $\pi'$ lisse et $S$ connexe. De plus, la réduction de Stein nous ramène au cas d'un morphisme à fibres connexes. Ainsi, on peut supposer le morphisme $\pi^{'}$ lisse et à fibres connexes sur une base $S$ normale et connexe. Dans ce cas,  les faisceaux  
${\rm I}\!{\rm R}^{2n}\pi_{*}{\Bbb C}$ et 
${\rm I}\!{\rm R}^{n}\pi_{*}\Omega^{n}_{Y'/Y}$ sont respectivement localement 
constant et localement libre de
rang 1 sur $Y$, de plus on a les isomorphismes canoniques, définis par 
intégration\par
\centerline{ ${\rm I}\!{\rm R}^{2n}\pi'_{*}{\Bbb C}\simeq {\Bbb C} $}\par
\centerline{ ${\rm I}\!{\rm R}^{n}\pi_{*}\Omega^{n}_{Y'/Y}\simeq{\mathcal 
O}_{Y},$}\smallskip
\noindent
  Par conséquent, pour tout point $y_{0}$ fixé de $Y$,  tout   voisinage ouvert, $V$,   de $\pi'^{-1}(y_{0})$ (que l'on peut 
choisir  $n$- complet d'après \cite{B4}) et toute forme  $Y$-relative, $\alpha$,  
$\bar\partial$-fermée de type $(n,n)$ 
sur $V$, il existe une fonction holomorphe
$f(y)$ définie sur un voisinage ouvert de $y_{0}$ dans $Y$ et $\beta$ 
une $(n,n-1)$ forme 
relative telle que
$$f(y).\omega_{1}^{n} - \alpha = \bar\partial \beta.$$
\noindent
 Choisissons $\alpha = \omega_{2}^{n}$ et $f$ telle que  
$f(y_{0})=1$, ce qui est possible puisque $\displaystyle
 \int_{Y'/Y}\omega_{1}^{n} = \int_{Y'/Y}\omega_{2}^{n} = 1$ et donc
$({\omega_{1}^{n}} - {\omega_{2}^{n}})\mid{\pi^{-1}(y)}$ est
 $\bar\partial$ exacte pour tout point $y$ de $Y$.  Alors, si $r$ est  la dimension de $Y$, on a, pour toute  forme ${\mathcal C}^{\infty}$ de type $(r-q,r)$ sur $Y$,  $\phi$,
$$\int_{Y'}\pi^{*}(\xi) \wedge {{\omega}_{2}}^{n}\wedge 
\pi'^{*}(\phi) =
\int_{Y'}\pi^{*}(\xi) \wedge \{f(y){{\omega}_{1}^{n}}- 
\bar\partial \beta\}
 \wedge \pi'^{*}(\phi)$$
\noindent Mais pour des raisons de $\bar\partial$ -fermeture et de type, on a
$$\int_{Y'}\pi^{*}(\xi) \wedge \bar\partial{ \beta\wedge 
\pi'^{*}(\phi)} = 
\int_{Y'}\bar\partial\{\pi^{*}(\xi) \wedge \beta\wedge 
\pi'^{*}(\phi)\}= \int_{Y'}d\{\pi^{*}(\xi) \wedge \beta\wedge \pi'^{*}(\phi)\}=0$$
\noindent Et donc
$$\int_{Y'}\pi^{*}(\xi) \wedge {{\omega}_{2}^{n}}\wedge 
\pi'^{*}(\phi) =
 \int_{Y'}\pi^{*}(\xi) \wedge f(y){{\omega}_{1}^{n}}\wedge 
\pi'^{*}(\phi) $$ ou
$${\pi'_{*}}(\pi^{*}(\xi)\mid_{Y'} \wedge {{\omega}_{2}^{n}}\mid_{Y'})={\pi'_{*}}(\pi^{*}(\xi)\mid_{Y'} \wedge f(y){{\omega}_{1}^{n}}\mid_{Y'})$$
Mais pour $\xi$ forme holomorphe, cette relation devient (cf \propositionref{P'2}) en raison de la normalisation choisie pour les formes de Kähler
$$i^{*}(\xi)=f(y)i^{*}(\xi)$$
D'où $f(y)=1,\,\,\forall\,y\in Y$ et donc 
$$ {\mathcal R}^{\pi,\omega_{1}}_{Z,Y}= {\mathcal R}^{\pi,\omega_{2}}_{Z,Y}$$
\indent {\bf (b) l'indépendance vis-à-vis du morphisme choisi.}\rm\par
\noindent
 Soient $\pi_{1}: Z_{1}\rightarrow Z$ (resp. $\pi_{2}: Z_{2}\rightarrow Z$) un morphisme k\"ahlérien  (resp. localement  k\"ahlérien) avec $Z_{1}$ (resp. $Z_{2}$) variété analytique complexe connexe  et $\omega_{1}$ (resp. $\omega_{2}$)  une forme de K\"ahler relative au morphisme $\pi_{1}$ (resp. $\pi_{2}$). On notera  $Y_{1}:=\pi^{-1 }_{1}(Y)$ (resp. $Y_{2}:=\pi^{-1 }_{2}(Y)$ ) la préimage totale de $Y$ relativement à $\pi_{1}$ (resp. $\pi_{2}$).  On veut montrer l'égalité 
$$ {\mathcal R}^{\pi_{1},\omega_{1}}_{Z_{1},Y_{1}}= {\mathcal R}^{\pi_{2},\omega_{2}}_{Z_{2},Y_{2}}$$
Pour cela, on remarque que la catégorie des morphismes k\"ahlériens  (resp. localement  k\"ahlériens) possèdent les propriétés requises pour établir cette invariance puisque cette  notion est stable par changement de base arbitraire, par composition et qu'une classe importante de tels morphismes est donnée par les projectifs (resp. localement projectifs) et, donc, en particulier par les modifications propres (cf [Bing]).  Ces propriétés nous permettent de construire, à partir de deux $Z$-morphismes un troisième $Z$-morphisme  dans cette catégorie puisqu'il nous suffit  de  considérer la désingularisée de l'espace réduit sous jacent au produit fibré $Z_{1}\times_{Z} Z_{2}$ . On obtient, alors, un diagramme commutatif
$$
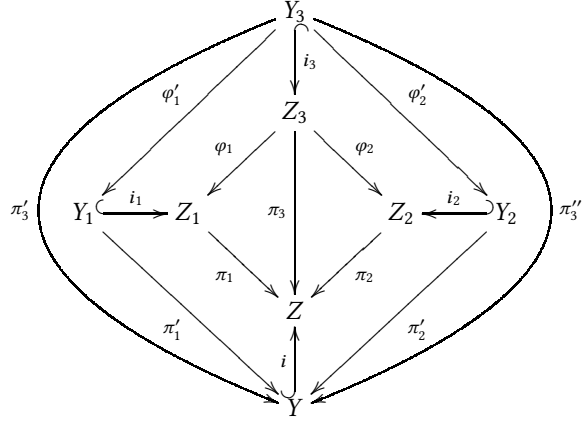
\xymatrix{&&Y_{3}\ar@/_8pc/[dddd]_{\pi'_{3}}\ar@/^8pc/[dddd]^{\pi''_{3}} \ar@{^{(}->}[d]^{i_{3}}\ar[lldd]_{\phi'_{1}}\ar[rrdd]^{\phi'_{2}}&&\\
&&Z_{3}\ar[dd]_{\pi_{3}}\ar[ld]_{\phi_{1}}\ar[rd]^{\phi_{2}}&&\\
Y_{1}\ar@{^{(}->}[r]^{i_{1}}\ar[rrdd]_{\pi'_{1}}&Z_{1}\ar[rd]_{\pi_{1}}&&Z_{2}\ar[ld]^{\pi_{2}}&Y_{2}\ar[lldd]^{\pi'_{2}}\ar@{_{(}->}[l]_{i_{2}}\\
&&Z&&\\
&&Y\ar@{^{(}->}[u]^{i}&&}$$
dans lequel  tous les morphismes du carré intérieur sont localement k\"ahlériens,   $Y_{3}$ est  naturellement construit à partir de l'espace réduit sous jacent au produit fibré $Y_{1}\times_{Y} Y_{2}$, $\phi'_{1}$ (resp.  $\phi'_{2}$) a même dimension générique que $\pi'_{2}$ (resp. $\pi'_{1}$). On note $n_{1}$ (resp. $n_{2}$) la dimension générique relative de  $\pi'_{1}$ (resp.  $\pi'_{2}$).\par\noindent
Notre but est de montrer l'égalité au sens des courants  $${\pi'_{1}}_{*}({\pi_{1}}^{*}(\xi)\wedge \omega_{1}^{n_{1}}|_{Y_{1}})={\pi'_{2}}_{*}({\pi_{2}}^{*}(\xi)\wedge \omega_{2}^{n_{2}}|_{Y_{2}})$$
qui, grâce à la commutativité du diagramme ci-dessus, revient à montrer
$${\pi'_{1}}_{*}({\pi_{1}}^{*}(\xi)\wedge \omega_{1}^{n_{1}}|_{Y_{1}})={\pi''_{3}}_{*}({\pi_{3}}^{*}(\xi)\wedge \omega_{3}^{n_{3}}|_{Y_{3}}),\,\,\,n_{3}=n_{1}+ n_{2}$$
On se ramène ainsi  à un diagramme du type
$$
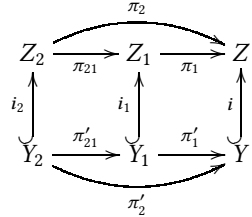
\xymatrix{Z_{2}\ar@/^1pc/[rr]^{\pi_{2}}\ar[r]_{\pi_{21}}&Z_{1}
\ar[r]_{\pi_{1}}&Z\\
Y_{2}\ar@{^{(}->}[u]^{i_{2}}\ar[r]^{\pi'_{21}}\ar@/_1pc/[rr]_{\pi'_{2}} &Y_{1}\ar@{^{(}->}[u]^{i_{1}}\ar[r]^{\pi'_{1}}&Y\ar@{^{(}->}[u]^{i}}$$
pour lequel, on va montrer l'égalité au sens des courants  $${\pi'_{1}}_{*}({\pi_{1}}^{*}(\xi)\wedge \omega_{1}^{n_{1}}|_{Y_{1}})={\pi'_{2}}_{*}({\pi_{2}}^{*}(\xi)\wedge \omega_{2}^{n_{2}}|_{Y_{2}})$$
Or, on a
$${\pi'_{2}}_{*}({\pi_{2}}^{*}(\xi)\mid_{Y_2}\wedge \omega_{2}^{n_{2}}|_{Y_{2}})={\pi'_{1}}_{*}{\pi'_{21}}_{*}({\pi'_{21}}^{*}({\pi_{1}}^{*}(\xi|_{Y_{1}})\wedge  \omega_{2}^{n_{2}}|_{Y_{2}})={\pi'_{1}}_{*}({\pi_{1}}^{*}(\xi|_{Y_{1}})\wedge{\pi'_{21}}_{*}( \omega_{2}^{n_{2}}|_{Y_{2}}))$$
D'après {\bf(a)}, la construction étant indépendante de la forme de K\"ahler choisie, on peut donc prendre $\omega_{2}$ telle qu'au niveau des classes de K\"ahler on ait
$\omega_{2}^{n_{2}}={\pi_{21}}^{*}(\omega_{1}^{n_{1}})\wedge \omega_{21}^{n_{21}}$,  $\omega_{21}$ une forme de K\"ahler relative au morphisme $\pi_{21}$. Dans ce cas,  grâce à la normalisation usuelle  ${\pi'_{21}}_{*}(\omega_{21}^{n_{21}})=1$,  on aboutit à l'égalité désirée\footnote{A noter que, relativement au grand diagramme initial,  la forme $\omega_{3}$ peut être définie au moyen de ${\phi_{1}}^{*}( \omega_{1})$ et  ${\phi_{2}}^{*}( \omega_{2})$ par construction de $Z_{3}$).}.\par\noindent
On peut remarquer que, comme on peut supposer tous les morphismes du carré extérieur du diagramme précédent géométriquement plats,  le théorème de Varouchas (cf  \cite{Va}) garantissant que  l'image directe par un morphisme géométriquement plat  d'une forme de K\"ahler est encore 
une forme de K\"ahler  montre que le courant ${\phi'_{1}}_{*}( \omega_{3}^{n_{3}}|_{Y_{3}})$ est la puissance $n_{1}$-ème d'une forme de K\"ahler.  \smallskip\noindent
Cette invariance vis-à-vis du morphisme choisi permet de prendre pour morphisme $\pi$ tout  morphisme de désingularisation ce que l'on supposera par la suite.\vskip 1mm\noindent
{\bf 3.4.3. Propriétés fonctorielles de cette restriction.}\vspace{1mm}

\noindent
{\bf (i) Compatibilité avec le cas lisse et avec la restriction usuelle sur les formes holomorphes.}\rm\par
\noindent 
 C'est l'objet de la  \propositionref{P'1} et  de la \propositionref{P'2}. \vspace{1mm}

\noindent 
{\bf (ii) Formule de transitivité.}\vspace{1mm}

\noindent
Il est clair,  d'après la construction générale, que le cas fondamental est celui où $T'$ est contenu dans le lieu singulier de $T$ et $T$ contenu dans le lieu singulier de $Z$.  Notre problème étant de nature local, on peut supposer $Z$, $T$ et $T'$ des ouverts de Stein et $\Sigma$ irréductible.\vspace{1mm}

\noindent
Il importe de regarder l'incidence de $T$ et du lieu singulier ${\rm Sing}(Z)$ de $Z$ que l'on classe selon le schéma:\vspace{1mm}

\indent{\bf {Cas 1: $T$ est contenu dans la partie régulière de $Z$.} }\vspace{1mm}

\noindent
On se ramène aux restrictions usuelles des formes holomorphes comme expliqué dans  la \propositionref{P'1}.\vspace{1mm}

\indent{\bf {Cas 2: $T$ n'est pas entièrement contenu dans ${\rm Sing}(Z)$  et  $T':=T\cap {\rm Sing}(Z)$.}} \vspace{1mm}

\noindent
On considère le diagramme
$$\xymatrix{\Sigma'\ar[r]^{\!\!\!\!\!\!\!\!\!{\sigma}'}\ar[d]_{\pi'} &\Sigma\ar[d]_{\pi}\ar[r]^{{\sigma}}&\tilde{Z}\ar[d]^{\phi\,\,\,\,\,\,\,\,\,\,\,\,\,\,\,(\spadesuit)}\\
T'\ar[r]^{ i'}&T\ar[r]^{ i}&Z}$$
avec $\pi$ modification propre, $\Sigma$ (resp.  $\Sigma'$) la préimage stricte (resp. totale)  de $T$ (resp. $T'$). A noter que l'on peut se ramener  au cas où $\Sigma$ est lisse (quitte à le désingulariser) avec $\sigma$ et $\sigma'$ des modifications propres.\vspace{1mm}

\noindent 
Alors, $${\mathcal R}_{Z,T}(\xi):=\xi_{T}:={\pi}_{*}( \phi^{*}(\xi)\!\mid_{\Sigma})$$
$${\mathcal R}_{Z,T'}(\xi):=\xi_{T'}:= \pi'_{*}( \phi^{*}(\xi)\!\mid_{\Sigma'}\wedge \omega^{n}\!\!\mid_{\Sigma'})$$
$${\mathcal R}_{T,T'}(\xi_T):={\xi}_{T,T'}:={\pi'_{*}}( {\pi}^{*}(\xi_{T})\!\!\mid_{\Sigma'}\wedge {\omega}^{n}\!\!\mid_{\Sigma'})$$
La  relation de transitivité que l'on cherche à établir est $${\mathcal R}_{Z,T'}={\mathcal R}_{T,T'}\circ {\mathcal R}_{Z,T}$$
Comme $\Sigma\setminus \Sigma'\simeq T\setminus T'$ et $\tilde{Z}\setminus E\simeq Z\setminus{\rm Sing}(Z)$,  la forme holomorphe  $\phi^{*}(\xi)\!\mid_{\Sigma})$ définit, par construction,  une section du faisceau $\Omega^{q}_{\Sigma}/{\rm Tors} $ (puisque non triviale sur un ouvert dense). Cela nous permet de supposer $\Omega^{q}_{\Sigma}$ sans torsion et, donc, le voir comme un sous faisceau de ${\mathcal L}^{q}_{T}$. De toute évidence, $\xi_{T}$ est une section du faisceau  $\pi_{*}(\Omega^{q}_{\Sigma})$. La factorisation naturelle (${\mathcal T}$ étant le sous faisceau de torsion)
$$\xymatrix{0\ar[r]&{\mathcal T}\ar[r]&\pi^{*}(\pi_{*}(\Omega^{q}_{\Sigma}))\ar[r]\ar[rd]&\pi^{*}(\pi_{*}(\Omega^{q}_{\Sigma}))/{\mathcal T}\ar[r]\ar@{^{(}->}[d]&0\\
&&&\Omega^{q}_{\Sigma}&}$$
montre que  $\pi^{*}(\xi_{T})$ ne peut pas être de torsion dans $\pi^{*}(\pi_{*}(\Omega^{q}_{\Sigma}))$ car sinon son  image dans $\Omega^{q}_{\Sigma}$ serait trivialement nulle!
Pour terminer avec les remarques évidentes, notons que l'on a les isomorphismes
$$\pi_{*}\pi^{*}(\pi_{*}(\Omega^{q}_{\Sigma}))\simeq\pi_{*}(\Omega^{q}_{\Sigma}),\,\, \pi_{*}\pi^{*}({\mathcal L}^{q}_{T})\simeq {\mathcal L}^{q}_{T}$$
\vspace{1mm} 

\noindent
puisque, par exemple, la factorisation (analogue de celle ci-dessus) 
$$\xymatrix{0\ar[r]&{\mathcal T}\ar[r]&\pi^{*}\pi_{*}({\mathcal L}^{q}_{\Sigma})\ar[r]\ar[rd]_{u}&\pi^{*}\pi_{*}({\mathcal L}^{q}_{\Sigma})/{\mathcal T}\ar[r]\ar@{^{(}->}[d]&0\\
&&&{\mathcal L}^{q}_{\Sigma}&}$$
montre que ${\rm Ker}u={\mathcal T}$ et, par conséquent, le morphisme induit
$$\pi_{*}\pi^{*}({\mathcal L}^{q}_{T})\rightarrow {\mathcal L}^{q}_{T}$$
est injectif. Or, on a a déja un morphisme naturel
$${\mathcal L}^{q}_{T}\rightarrow \pi_{*}\pi^{*}({\mathcal L}^{q}_{T})$$
nécessairement injectif puisque c'est un isomorphisme générique et ${\mathcal L}^{q}_{T}$ est sans torsion .  D'où l'isomorphisme
$$\pi_{*}\pi^{*}({\mathcal L}^{q}_{T})\simeq {\mathcal L}^{q}_{T}$$
On montre de la même façon  le premier isomorphisme.\vspace{1mm}

\noindent 
Ceci étant dit, on peut  écrire ${\pi}^{*}(\xi_{T})=(\phi^{*}(\xi))\!\!\!\mid_{\Sigma} +\chi$ où $\chi$ une $q$-forme holomorphe à support dans $\Sigma'$, sans torsion sur $\Sigma'$ (et  qui rencontre le diviseur exceptionnel). \vspace{1mm}

\noindent
Alors, si $\Sigma'_{0}$ est  le lieu régulier de $\pi'$ et quitte  à rétrécir  $T'_{0}$, on peut supposer $\Sigma'_{0}\simeq T'_{0}\times V$. Dans ce cas, le théorème d'intégration de Lelong nous donne
$$\int_{\Sigma'}\psi\wedge \chi\wedge\omega^{n}\!\!\mid_{\Sigma'}=\int_{\Sigma'_{0}}\pi'^{*}(\psi)\wedge \chi\wedge\omega^{n}\!\!\mid_{\Sigma'_{0}}$$
pour une fonction test $\psi$ de type $(r-q, r)$ sur $T'$. Mais cela montre clairement que $\chi\wedge \omega^{n}\!\!\mid_{\Sigma'_{0}} $ est non triviale si et seulement si  $\chi$ s'exprimer qu'en fonction des coordonnées de $T'$ c'est-à-dire qu'il existe, $\chi'_{0}$, une $q$-forme holomorphe sur $T'_{0}$ telle que 
$\chi=\pi'^{*}(\chi'_{0})$. En considérant une extension, $\chi'$, quelconque de $\chi'_{0}$ à $T'$, on obtient la relation
$${\pi}^{*}(\xi_{T})\!\!\mid_{\Sigma'_{0}}=(\phi^{*}(\xi))\!\!\mid_{\Sigma'_{0}} +\pi'^{*}(\chi')\!\!\mid_{\Sigma'_{0}}$$
En supposant, les faisceaux de formes holomorphes $\Omega^{\bullet}_{\Sigma'}$ et $\Omega^{\bullet}_{T'}$ sans torsion, cette relation ayant lieu sur la partie dense 
$\Sigma'_{0}$ s'étend par absence de torsion à $\Sigma'$ tout entier. Ainsi,
$${\pi}^{*}(\xi_{T})\!\!\mid_{\Sigma'}=(\phi^{*}(\xi))\!\!\mid_{\Sigma'} +\pi'^{*}(\chi')$$
qui s'écrit aussi,  pour une extension $\chi"$ de $\chi'$ à $T$,
$${\pi}^{*}(\xi_{T})\!\!\mid_{\Sigma'}=(\phi^{*}(\xi))\!\!\mid_{\Sigma'} +\pi^{*}(\chi")\!\!\mid_{\Sigma'}$$
et, donc, 
$$({\pi}^{*}(\xi_{T} -\chi")-\phi^{*}(\xi))\!\!\mid_{\Sigma})\!\!\mid_{\Sigma'}=0$$
Mais comme cette relation est vérifiée sur $\Sigma\setminus \Sigma'$, on en déduit
$${\pi}^{*}(\xi_{T} -\chi")=\phi^{*}(\xi)\!\!\mid_{\Sigma}$$
ou de façon équivalente (on peut supposer $T$ normal mais ce n'est pas nécessaire) 
$$\xi_{T} -\chi"=\pi_{*}(\phi^{*}(\xi))\!\!\mid_{\Sigma}:=\xi_{T}$$
qui donne $\chi"=0$ et, par suite, $$\chi=0$$
D'où la relation $${\pi}^{*}(\xi_{T})=\phi^{*}(\xi)\!\!\mid_{\Sigma}$$
et donc $$\xi_{T'}=\xi_{T,T'}$$
qui n'est autre que la relation de transivité désirée.\vspace{1mm}

\indent {\bf{Cas 3: $T$ contenu dans ${\rm Sing}(Z)$.}} \vspace{1mm}

\noindent
On considère une désingularisation (ou tout autre modification propre) $\phi:\tilde{Z}\rightarrow Z$ (resp. $\nu:\tilde{T}\rightarrow T$) de $Z$ (resp. $T$) et le diagramme  
$$\xymatrix{\tilde{\Sigma}'\ar[rr]^{{\tilde \sigma}}\ar[dd]_{\tilde{\pi'}}\ar[rd]_{\Theta'}&&
\tilde{\Sigma}\ar[rd]^{\Theta}\ar[dd]^{\tilde\pi}&&\\
&\Sigma'\ar[rr]^{\!\!\!\!\!\!\!\!\!{\sigma}'}\ar[dd]_{\pi'} &&\Sigma\ar[dd]_{\pi}\ar[r]^{{\sigma}}&\tilde{Z}\ar[dd]^{\phi\,\,\,\,\,\,\,\,\,\,\,\,\,\,\,(\spadesuit)}\\
 \tilde{T}'\ar[rr]\ar[rd]_{\nu'}&&\tilde{T}\ar[rd]^{\nu}&&\\
&T'\ar[rr]^{{ i}'}&&T\ar[r]^{{ i}}&Z}$$
dans lequel  les morphismes horizontaux sont des immersions fermées, $\Sigma$ (resp. $\Sigma'$) la  préimage totale de $T$ (resp. $T'$) pour la modification $\phi$ et $\tilde{T}'$ celle de $T'$ pour $\nu$.\par\noindent  Les morphismes étant localement k\"ahlérien, on considère des formes de K\"ahler  $\omega$, $\tilde\omega$ et  $\hat\omega$ relative aux morphismes $\phi$, $\nu$ et $\pi$ respectivement.
Alors, comme tous les  morphismes  de ces diagrammes sont  propres, les images directes au sens des courants
$${\mathcal R}_{Z,T}(\xi):=\xi_{T}:={\pi}_{*}( \phi^{*}(\xi)\!\mid_{\Sigma}\wedge {\omega}^{n}\!\!\mid_{\Sigma})$$
$${\mathcal R}_{Z,T'}(\xi):=\xi_{T'}:= \pi'_{*}( \phi^{*}(\xi)\!\mid_{\Sigma'}\wedge \omega^{n}\!\!\mid_{\Sigma'})$$
$${\mathcal R}_{T,T'}(\xi_T):={\xi}_{T,T'}:={\nu}'_{*}( {\nu}^{*}(\xi_{T})\!\!\mid_{{\tilde T}'}\wedge {\tilde\omega}^{m}\!\!\mid_{{\tilde T}'})$$
$$\widetilde{\xi}_{T,T'}:={\pi'_{*}}( \pi^{*}(\xi_{T})\!\mid_{\Sigma'}\wedge \omega^{n}\!\!\mid_{\Sigma'})$$
définissent des courants de type $(q,0)$, $\bar\partial$-fermés induisant,  modulo torsion, des sections des faisceaux $\omega^{q}_{T}$ et  $\omega^{q}_{T'}$. Avec ces notations, la formule de transitivité recherchée se traduit par les égalités 
$${\mathcal R}_{Z,T'}={\mathcal R}_{T,T'}\circ {\mathcal R}_{Z,T}$$
On a aussi
$$\xi_{T,T'}=\widetilde{\xi}_{T,T'}$$
Pour les prouver, on commence par traiter le cas où:\par\indent
{\bf(a) $\pi:\Sigma\rightarrow T$ et $\nu': \tilde{T'}\rightarrow \tilde{T}$ géométriquement plats à fibres de dimension $n$ et $m$ respectivement.}\vspace{2mm}

\noindent
La notion de platitude géométrique étant  stable par changement de base, les morphismes verticaux $\pi'$, $\tilde\pi$, $\tilde\pi'$  sont géométriquement plats à fibres de même dimension $n$ puisque $\pi$. Il en est de même de $\Theta'$ en raison de la platitude géométrique de $\nu$. Dans ce contexte, les images directes définies au sens des courants  coincident avec l'intégration sur les fibres, opérations dont la formation commute aux changements de base.\par\noindent
Comme le diagramme précédent est commutatif et que  $\pi'$ se déduit de $\pi$ dans le changement de base $T'\rightarrow T$ (par définition $\Sigma'=T'\times_{Z}\tilde{Z}=T'\times_{T}\Sigma$ et  $\Sigma'=T\times_{Z}\tilde{Z}$),  on a, pour toute forme $\alpha$ sur $\Sigma$:
$$i'^{*}(\int_{\Sigma/T}(\alpha\wedge\omega^{n}\mid_{\Sigma} ))=\int_{\Sigma'/T'}{\sigma'}^{*}(\alpha\wedge\omega^{n}\mid_{\Sigma})=\int_{\Sigma'/T'}(\alpha\wedge\omega^{n})\mid_{\Sigma'}$$
qui, pour $\alpha:= \phi^{*}(\xi)\mid_{\Sigma}$, s'écrit aussi
$$i'^{*}(\xi_{T})=\pi'_{*}(\phi^{*}(\xi)\mid_{\Sigma'}\wedge\omega^{n}\mid_{\Sigma'})=\xi_{T'}:={\mathcal R}_{Z,T'}(\xi)$$
c'est-à-dire $i'^{*}\circ{\mathcal R}_{Z,T}={\mathcal R}_{Z,T'}$.  En fait, dans ce contexte, on a   $i'^{*}={\mathcal R}_{T,T'}$ qui est essentiellement dû à la propriété de commutation aux changements de base arbitraires de l'intégration sur une famille analytique de cycles.  En effet, pour montrer l'égalité 
$$\pi'_{*}( \phi^{*}(\xi)\!\!\mid_{\Sigma'}\wedge \omega^{n}\!\!\mid_{\Sigma'})= {\nu}'_{*}( {\nu}^{*}(\xi_{T})\!\!\mid_{{\tilde T}'}\wedge {\tilde\omega}^{m}\!\!\mid_{{\tilde T}'})$$
on utilise cette propriété pour écrire 
$$\nu^{*}(\xi_{T}):=\nu^{*}(\int_{\Sigma/T}(\phi^{*}(\xi)\wedge\omega^{n})\!\mid_{\Sigma} )=\int_{\tilde\Sigma/\tilde T}{\Theta}^{*}((\phi^{*}(\xi)\wedge\omega^{n})\!\mid_{\Sigma})$$
$$\nu^{*}(\xi_{T})\!\!\mid_{{\tilde T}'}:=\int_{\tilde\Sigma'/\tilde T'}{\Theta}'^{*}((\phi^{*}(\xi)\wedge\omega^{n})\!\mid_{\Sigma'})$$
ou 
$${\nu}^{*}(\xi_{T})={\nu}^{*}( \pi_{*}((\phi^{*}(\xi)\mid_{\Sigma}\wedge \omega^{n}\!\!\mid_{\Sigma}))={\tilde\pi}_{*}(\Theta^{*}((\phi^{*}(\xi)\mid_{\Sigma}\wedge \omega^{n}\!\!\mid_{\Sigma}))$$
$${\nu}^{*}(\xi_{T})\!\!\mid_{{\tilde T}'}={\tilde\pi}'_{*}((\Theta^{*}((\phi^{*}(\xi)\!\!\mid_{\Sigma}\wedge \omega^{n}\!\!\mid_{\Sigma}))\!\!\mid_{{\tilde\Sigma}'})={\tilde\pi}'_{*}( \Theta'^{*}(( \phi^{*}(\xi)\!\!\mid_{\Sigma'}\wedge \omega^{n}\!\!\mid_{\Sigma'}) ))$$
Alors, la formule de projection pour un morphisme propre permet d'écrire
$${\nu}'_{*}( {\nu}^{*}(\xi_{T})\!\!\mid_{{\tilde T}'}\wedge {\tilde\omega}^{m}\!\!\mid_{{\tilde T}'})={\nu}'_{*}({\tilde\pi}'_{*}( \Theta'^{*}((\phi^{*}(\xi)\!\!\mid_{\Sigma'}\wedge \omega^{n}\!\!\mid_{\Sigma'}) )\wedge {\tilde\omega}^{m}\!\!\mid_{{\tilde T}'})$$
$$= \pi'_{*}\Theta'_{*}\{\Theta'^{*}((\phi^{*}(\xi)\!\!\mid_{\Sigma'}\wedge \omega^{n}\!\!\mid_{\Sigma'} )\wedge {\tilde\pi}'^{*}({\tilde\omega}^{m}\!\!\mid_{{\tilde T}'})\}= \pi'_{*}\{( (\phi^{*}(\xi)\!\!\mid_{\Sigma'}\wedge \omega^{n}\!\!\mid_{\Sigma'})\wedge\Theta'_{*}({\tilde\pi}'^{*}({\tilde\omega}^{m}\!\!\mid_{{\tilde T}'})\}$$
 La formule du changement de base pour l'intégration le long des fibres de $\Theta'$ nous donne 
 $${\Theta'_{*}{\tilde\pi}'^{*} ({\tilde\omega}^{m}\!\!\mid_{{\tilde T}'})=\pi'^{*}({\nu}'_{*}{\tilde\omega}^{m}\!\!\mid_{{\tilde T}'})}$$
 et donc, en vertu de la formule de projection, $$\nu'_{*}( {\nu}^{*}(\xi_{T})\!\!\mid_{{\tilde T}'}\wedge {\tilde\omega}^{m}\!\!\mid_{{\tilde T}'})=\pi'_{*}(( \phi^{*}(\xi)\!\!\mid_{\Sigma'}\wedge \omega^{n}\!\!\mid_{\Sigma'})\wedge{\nu}'_{*}{\tilde\omega}^{m}\!\!\mid_{{\tilde T}'}$$
Alors, en normalisant comme de coutume la forme de K\"ahler par  ${\nu}'_{*}({\tilde\omega}^{m}\!\!\mid_{{\tilde T}'})=1$, on obtient
$${\mathcal R}_{T,T'}(\xi_T)={\xi}_{T'}$$
qui est exactement la formule de transitivité cherchée. 
\vspace{1mm}

\noindent
Comme notre  construction de la restriction ne dépend pas de la forme de K\"ahler choisie,  on en en déduit la relation de transitivité en toute généralité et 
$$i'^{*}(\xi_{T})=\xi_{T'}={\tilde\xi}_{T,T'}$$
{\bf(b)} $\pi\!\!:\!\!\Sigma\rightarrow T$ {\bf et}  $\nu'\!\!:\!\!\tilde{T'}\rightarrow \tilde{T}$ {\bf{propres surjectifs génériquement $n$ et $m$-équidimensionnels respectivement:}}\vspace{2mm}

\noindent
Pour nous ramener au cas géométriquement plat, observons que l'on peut  supposer:\vspace{1mm}

\indent 
$\bullet$ {\it  $T'$  est le lieu de non platitude géométrique (resp. de non platitude algébrique) de $\pi$}: \vspace{1mm}

\noindent
 Comme $\pi$ est un morphisme propre d'espaces complexes réduits, on peut commencer par le platifier grâce au théorème de Hironaka (\cite{H1}). Dans le diagramme ($\spadesuit$), $\nu$ est la modification sur la base $T$ pour transformer $\pi$ en le morphisme plat $\tilde{\pi}$. Alors, si:\vspace{1mm}
 
 \indent $\bullet$ {\it{$T'$ n'est pas contenu dans le lieu de la modification}}, il sera contenu dans l'ouvert dense du lieu de platitude et, donc, de platitude géométrique. Dans ce cas, $\tilde{T'}$ est la préimage stricte et $\nu'$ est une modification propre.  On se retrouve dans la situation du cas {\bf(a)} où les morphismes sont géométriquement plats car on peut remplacer $\pi$ et $\pi'$ par les morphismes plats (et donc géométriquement plats) $\tilde\pi$ et $\tilde\pi'$ respectivement. \vspace{1mm}

\indent
 $\bullet$ {\it{$T'$ est contenu dans le lieu de la  modification}} alors $\nu'$ est propre génériquement $m$-équidimensionnel. Quitte à modifier $T'$, on se ramène au cas géométriquement plat et on termine comme dans {\bf(a)}.\vspace{1mm}

\noindent 
 De plus, on peut, si on le désire, supposer $\Sigma:= \pi^{-1}(T)$ lisse. En effet, s'il ne l'était pas,  on le désingularise grâce au théorème de Hironaka pour obtenir le diagramme commutatif
$$\xymatrix{{\tilde\Sigma}'\ar[r]^{\tilde{i}'}\ar[d]_{\psi'}&{\tilde{\Sigma}}\ar[d]_{\psi}&\\
 \Sigma'\ar[r]^{{\bf i}'}\ar[d]_{\pi'} &\Sigma\ar[r]^{\bf{i}}\ar[d]_{\pi}&{\tilde Z}\ar[d]_{\phi}\\
 T'\ar[r]_{{i'}}&T\ar[r]_{i}&Z }$$
Considérons, alors,  $\omega$ et $\tilde\omega$  deux formes de K\"ahler sur $\tilde{Z}$ et
$\tilde{\Sigma}$, relatives aux morphismes $\phi$ et ${\pi}o\psi$. Comme précédemment, on désigne par $n$ la dimension  générique  relative du morphisme $\pi$, on pose $ \tilde{\pi}':={\pi'}o\psi'$ et $\tilde{\pi}:={\pi}o\psi$. Comme la construction de notre restriction est indépendante de la classe de K\"ahler choisie, on peut  prendre
$[\tilde\omega] = [\psi^{*}(\sigma^{*}\omega)]$. Il est facile de voir que les définitions  des formes $\xi_{T}$,  $\xi_{T'}$ et ${\tilde\xi}_{T'}$  relatives à cette configuration sont les mêmes que celles données plus haut. En effet,  un calcul simple en termes d'images directes au sens des courants (utilisant la formule de projection pour un morphisme propre)  montre, par exemple, que 
$$\xi_{T}:= \tilde{\pi}_{*}( \psi^{*}(\phi^{*}(\xi)\mid_{\Sigma})\wedge \tilde{\omega}^{n}\mid_{\tilde\Sigma})= {\pi}_{*}( \phi^{*}(\xi)\mid_{\Sigma}\wedge {\omega}^{n}\mid_{\Sigma})$$
puisque
$$\langle{(\pi o \psi)_{*}( \psi^{*}(\phi^{*}(\xi)\mid_{\Sigma}\wedge \tilde{\omega}^{n}\mid_{\tilde\Sigma})), \beta}\rangle= \langle{ \psi^{*}(\phi^{*}(\xi)\mid_{\Sigma}\wedge \tilde{\omega}^{n}\mid_{\tilde\Sigma})), \psi^{*}(\pi^{*}\beta)}\rangle$$
$$\indent\indent\indent\indent\indent\indent\indent\indent\indent\indent\indent=\langle{\psi_{*}[\tilde\Sigma],\phi^{*}(\xi)\mid_{\Sigma}\wedge \tilde{\omega}^{n}\mid_{\tilde\Sigma}))\wedge\pi^{*}\beta) }\rangle$$
et que  $\psi_{*}([\tilde\Sigma])=[\Sigma]$.\vspace{1mm}

\noindent
Prouvons, enfin, les égalités annoncées en les déduisant du cas géométriquement plat. Signalons que  le morphisme $\nu$ du diagramme ($\spadesuit$)  correspond, maintenant,  à une suite finie d'éclatements de centres situés dans $T'$  réalisant  " l'applatissement  géométrique "  $\tilde\pi$ du morphisme $\pi$. Il  s'en suit que $\tilde\pi'$ est aussi géométriquement plat.  On peut supposer ${\tilde T}$ lisse sans enfreindre la généralité. \par\noindent
Alors, d'après le cas {\bf(i)}, on sait que   les images directes définies au sens des courants ou par intégration  le long des fibres
$$\xi_{\tilde T}:={\tilde\pi}_{*}( \Theta^{*}(\phi^{*}(\xi)\!\mid_{\Sigma}\wedge {\omega}^{n}\!\!\mid_{\Sigma}))$$
$$\xi_{\tilde T'}:= {\tilde\pi}'_{*}( \Theta^{*}(\phi^{*}(\xi)\!\mid_{\Sigma}\wedge {\omega}^{n}\!\!\mid_{\Sigma})\!\mid_{\tilde\Sigma'})$$
$$\tilde{\xi}_{T,T'}:= {\tilde\pi}'_{*}( {\tilde\pi}^{*}(\xi_{\tilde T})\!\mid_{{\tilde\Sigma'}}\wedge {\hat\omega}^{n}\!\!\mid_{{\tilde\Sigma'}})$$
vérifient  les égalités
$\xi_{\tilde T}\!\!\mid_{\tilde T'}=\xi_{\tilde T'}=\tilde{\xi}_{T,T'}$ ( $\xi_{\tilde T}\!\!\mid_{\tilde T'}$ est alors la restriction  usuelle des formes holomorphes) traduisant la relation
$${\mathcal R}_{Z, \tilde{T}'}={\mathcal R}_{\tilde{T}, \tilde{T}'}\circ {\mathcal R}_{Z, \tilde{T}}$$
On veut montrer  que cela entraîne  l'égalité des formes $\xi_{T'}$, $\tilde{\xi}_{T,T'}$ et $\xi_{T,T'}$ définies génériquement par les formules d'intégration ou globalement par image directe au sens des courants. Or cela découle simplement du fait que $\xi_{T}:=\nu_{*}\xi_{\tilde T}$ et $\xi_{T'}:=\nu'_{*}\xi_{\tilde T'}$ puisque, par exemple, 
$$\nu_{*}(\xi_{\tilde T})=\nu_{*}({\tilde\pi}_{*}( \Theta^{*}(\phi^{*}(\xi)\!\mid_{\Sigma}\wedge {\omega}^{n}\!\!\mid_{\Sigma}))= {\pi}_{*}\Theta_{*} (\Theta^{*}(\phi^{*}(\xi)\!\mid_{\Sigma}\wedge {\omega}^{n}\!\!\mid_{\Sigma})):=\xi_{T}$$
$\Theta$ étant une modification propre ($\Sigma$ peut-être supposé lisse si on le désire). Si on veut éviter l'hypothèse de lissité (ou de normalité) sur $\Sigma$, on utilise la formule de projection
$$\Theta_{*} (\Theta^{*}(\phi^{*}(\xi)\!\mid_{\Sigma}\wedge {\omega}^{n}\!\!\mid_{\Sigma}))=\phi^{*}(\xi)\!\mid_{\Sigma}\wedge \Theta_{*}\Theta^{*}({\omega}^{n}\!\!\mid_{\Sigma})$$
et préciser que $\Theta_{*}\Theta^{*}({\omega}^{n}\!\!\mid_{\Sigma})$ est encore une forme de Kähler relative  pour le morphisme $\pi$.\vspace{1mm}

\noindent De façon similaire, on a   
$$\nu'_{*}\xi_{\tilde T'}=\nu'_{*}{\tilde\pi'}_{*}( \Theta'^{*}((\phi^{*}(\xi)\!\mid_{\Sigma'})\wedge {\omega}^{n}\!\!\mid_{\Sigma'}))={\pi'_{*}}(\phi^{*}(\xi)\!\mid_{\Sigma'}\wedge \Theta_{*}{\omega}^{n}\!\!\mid_{\Sigma'})=\xi_{T'}$$
On en déduit, alors, que 
$${\mathcal R}_{{T}, {T}'}(\xi_{T})=\xi_{T'}$$
puisque $${\mathcal R}_{{T}, {T}'}(\xi_{T})={\mathcal R}_{{T}, {T}'}(\nu_{*}\xi_{\tilde T})=\nu'_{*}((\nu^{*}\nu_{*}\xi_{\tilde T})\!\mid_{\tilde T'}\wedge {\tilde\omega}^{m}\!\mid_{\tilde T'})$$
qui, d'après le théorème d'intégration de Lelong et ce qui précède, donne
$${\mathcal R}_{{T}, {T}'}(\xi_{T})=\nu'_{*}(\xi_{\tilde T})\!\mid_{\tilde T'}\wedge {\tilde\omega}^{m}\!\mid_{\tilde T'})=\nu'_{*}(\xi_{\tilde T'}\wedge {\tilde\omega}^{m}\!\mid_{\tilde T'})=\xi_{\tilde T'}$$
en choisissant la normalisation habituelle pour la forme de Kähler. \vspace{1mm}

\noindent On peut aussi remarquer que les courants ${\mathcal R}_{{T}, {T}'}(\xi_{T})$ et $\xi_{T'}$
définissent des sections globales du faisceau ${\mathcal L}^{q}_{T'}$. Comme elles  coincident génériquement
 et que ce faisceau est sans torsion, elles coincident globalement.\vspace{2mm}

\noindent 
De plus, on a aussi
$$\nu'_{*}\tilde{\xi}_{\tilde T, \tilde T'} =\tilde{\xi}_{T,T'}$$
En effet, comme (toujours au sens des courants comme il a été fait dans le cas géométriquement plat)) 
$$\tilde{\xi}_{\tilde T, \tilde T'}={\tilde\pi'}_{*}( {\tilde\pi}^{*}(\xi_{\tilde T})\!\mid_{{\tilde\Sigma'}}\wedge {\hat\omega}^{n}\!\!\mid_{{\tilde\Sigma'}})={\tilde\pi'}_{*}(\Theta^{*}(\phi^{*}(\xi))\!\mid_{{\tilde\Sigma'}}\wedge {\hat\omega}^{n}\!\!\mid_{{\tilde\Sigma'}}) $$
et donc
$$\tilde{\xi}_{\tilde T, \tilde T'}={\tilde\pi'}_{*}(\Theta'^{*}(\phi^{*}(\xi))\!\mid_{{\Sigma'}}\wedge {\hat\omega}^{n}\!\!\mid_{{\tilde\Sigma'}}) $$
D'où
$$\nu'_{*}\tilde{\xi}_{\tilde T, \tilde T'}={\pi'}_{*}( \phi^{*}(\xi))\!\mid_{{\Sigma'}}\wedge \Theta'_{*}{\hat\omega}^{n}\!\!\mid_{{\tilde\Sigma'}})= {\pi}'_{*}(\pi^{*}(\xi_{T})\!\mid_{{\tilde\Sigma'}}\wedge \Theta'_{*}{\hat\omega}^{n}\!\!\mid_{{\tilde\Sigma'}})=\tilde{\xi}_{T,T'}\,\,\blacksquare$$
\vspace{1mm}

\section{\color{blue}{La preuve du \theoremref{T1}.}}\vspace{1mm}

\noindent
  On veut montrer qu'à tout morphisme d'espaces complexes réduits  $f:Y\rightarrow Z$ est associé un unique morphisme de faisceaux
${\bf f}^{*}:{\mathcal L}^{\bullet}_{Z}\rightarrow f_{*}{\mathcal L}^{\bullet}_{Y} $
prolongeant l'image réciproque des formes usuelles et de construction compatible avec la composition de certains morphismes.\vspace{1mm}

\noindent 
Il est clair que la construction de cette image réciproque qui est de nature locale sur la source et la base se déduit directement du \theoremref{T4} en utilisant les factorisations locales 
via le graphe
$$\xymatrix{Y\ar@{^{(}->}[r]^{\sigma}\ar[rd]_{f}&Y\times Z\ar[d]^{p}\\
&Z}$$
On peut supposer les espaces  irréductibles  et les  morphismes seront  supposés surjectifs. \vspace{1mm}

\noindent
D'après le \theoremref{T4}, on a un morphisme de restriction ${\bf \sigma}^{*}:   {\mathcal L}^{\bullet}_{Y\times Z}\rightarrow \sigma_{*}{\mathcal L}^{\bullet}_{Y}$.  Pour la projection canonique, la construction de l'image réciproque ${\bf p}^{*}: {\mathcal L}^{\bullet}_{Y}\rightarrow p_{*}{\mathcal L}^{\bullet}_{Y\times Z}$ est évidente  et on vérifie facilement que $ {\bf \sigma}^{*}\circ{\bf p}^{*}={\bf f}^{*}$.\vspace{1mm}

\noindent
{\bf(ii)} On se ramène à regarder le diagramme 
$$
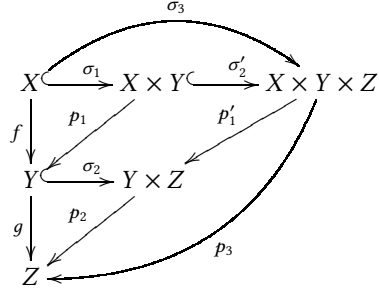
\xymatrix{X\ar[d]_{f}\ar@/^2pc/[rr]^{\sigma_3}\ar@{^(->}[r]^{\sigma_1}&X\times Y\ar[ld]_{p_1}\ar@{^(->}[r]^{\sigma'_2}&X\times Y\times Z\ar[ld]_{p'_1}\ar@/^2pc/[lldd]^{p_3}\\
Y\ar[d]_{g}\ar@{^(->}[r]^{\sigma_2}&Y\times Z\ar[ld]_{p_2}&\\
Z&&}$$
et vérifier l'égalité 
$$\bf{{\sigma}}^{*}_{1}({\bf p}_{1}^{*}( \bf{{\sigma}}^{*}_{2}({\bf p}_{2}^{*}\xi)))=\bf{{\sigma}}^{*}_{3}({\bf p}_{3}^{*}(\xi)) $$
qui résulte facilement du fait que, pour un diagramme cartésien
$$
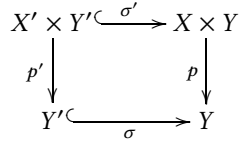
\xymatrix{X'\times Y'\ar@{^(->}[r]^{\sigma'}\ar[d]_{p'}&X\times Y\ar[d]_{p}\\
Y'\ar@{^(->}[r]_{\sigma}&Y}$$
où $\sigma$ (resp. $\sigma'$) et $p$ (resp. $p'$) désignent un plongement fermé et la projection canonique, on a
$${\bf \sigma'}^{*}\circ{\bf p}^{*}= {\bf p'}^{*}\circ {\bf \sigma}^{*}\,\blacksquare$$
\begin{rem} Cette image réciproque  semble souffrir d'une certaine incompatibilité avec la commutativité des diagrammes dans le sens où si 
$$
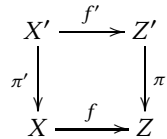
\xymatrix{{X'}\ar[d]_{\pi'}\ar[r]^{f'}&{Z'}\ar[d]^{\pi}\\
 X\ar[r]^{f}& Z}$$
est un tel diagramme d'espaces complexes réduits, il ne semble pas clair que 
$${\bf \pi'}^{*}\circ {\bf f}^{*}={\bf f'}^{*}\circ {\bf \pi}^{*}$$
La première impression est qu'une telle égalité ne peut avoir lieu que modulo certains noyaux d'opérateurs intégraux. Par exemple si $$
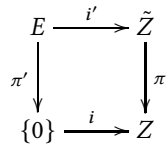
\xymatrix{{E}\ar[d]_{\pi'}\ar[r]^{i'}&\tilde{Z}\ar[d]^{\pi}\\
 \{0\}\ar[r]^{i}& Z}$$
 est l'éclatement en l'origine d'un espace complexe réduit donné de dimension $m$, toute forme de degré strictement supérieur à $1$ donne trivialement une restriction nulle. Mais, il se peut qu'elle donne par pull back sur $\tilde Z$, une forme dont la restriction à $E$ soit non triviale et dont l'intégrale sur $E$ du cup produit avec la puissance $m$-ème de la forme de Kähler relative soit nulle. On voit donc qu'il faut au moins quotienter par le noyau de cette opération d'intégration sur les fibres. 
\end{rem} 

\vspace{2mm}

\noindent
\phantomsection\addcontentsline{toc}{part}{La preuve du \theoremref{T2}.}
%
%
\svnid{$Id: S07-proof-setup.tex 269 2020-01-20 11:28:53Z kebekus $}

\section{\color{blue}{Image directe supérieure pour un morphisme géométriquement plat et ses conséquences.}}

\subversionInfo

\approvals{Mohamed & yes}
\par\vspace{2mm}

On voudrait associer à tout morphisme $\pi:X\rightarrow S$ possédant certaines propriétés de bonnes variances, au moins, au niveau des fibres, un morphisme de faisceaux de ${\mathcal O}_{S}$-modules (cohérents si $\pi$ est propre)
$${\mathcal T}^{k}_{\pi,{\mathcal L}}:{\rm I}\!{\rm R}^{n}\pi_{!}{\mathcal L}^{n+k}_{X}\rightarrow{\mathcal L}^{k}_{S} $$
compatible aux restrictions ouvertes sur $X$ et $S$, aux changements de base admissibles en un certain sens et aux morphismes ${\mathfrak T}^{k}_{\pi,\omega}:{\rm I}\!{\rm R}^{n}\pi_{!}{\omega}^{n+k}_{X}\rightarrow{\omega}^{k}_{S} $  du  {\it théorème 1} de \cite{K3}. Nous proposons, dans ce qui suit, deux constructions l'une plus géométrique et reposant sur une intégration le long des fibres et l'autre plus conceptuelle découlant de l'existence du morphisme trace défini au niveau des complexes dualisants (cf \cite{K3}). \vspace{1mm} 

\noindent La première approche consiste à la construire explicitement pour un morphisme géométriquement plat et en déduire, à moindre frais, la construction pour un morphisme propre surjectif arbitraire, puis, grâce aux théorèmes d'aplatissement locaux, le cas d'un morphisme surjectif génériquement ouvert. Pour cela, nous rappelons 
\Th{}{}{\color{blue}{Aplatissement local}}\label{Aplat2}\cite{H2}: \vspace{1mm}  \noindent Soit $\pi:X\rightarrow Y$ un morphisme d'espaces analytiques complexes avec $Y$ réduit. Soit $y\in Y$ et $L$ un sous ensemble compact de la fibre $\pi^{-1}(y)$. Alors, il existe un nombre fini de morphismes $\sigma_{\alpha}:Y_{\alpha}\rightarrow Y$ chacun d'entre eux composé d'un nombre fini d'éclatement locaux dont les centres sont situés dans des sous ensembles d'intérieur vides et tels que:\vspace{1mm}  
 
 \noindent
 {\bf(i)} Pour chaque $\alpha$, il existe un sous ensemble compact $K_{\alpha}$ de $Y_{\alpha}$ tel que $\displaystyle{\bigcup_{\alpha}\sigma_{\alpha}(K_{\alpha})}$ soit un voisinage compact de $y$ dans $Y$.\vspace{1mm} 
 
 \noindent
 {\bf(ii)} Pour chaque $\alpha$, la transformée stricte $\pi_{\alpha}:X_{\alpha}\rightarrow Y_{\alpha}$
 de $\pi$ par $\sigma_{\alpha}$ est plate aux points correspondant aux points de $L$.\rm
 \Th{}{}{\color{blue}{Aplatissement local géométrique}}\label{Aplat3}\cite{Si}: \vspace{1mm} 
\noindent Soient $\pi:X\rightarrow Y$ un morphisme génériquement ouvert d'espaces analytiques complexes  réduits et ${\rm K}$ un compact de $X$. Alors, il existe un espace complexe normal $Y'$, une application holomorphe $\sigma:Y'\rightarrow Y$ et un voisinage ouvert relativement compact $V$ de ${\rm K}$ tels que:\vspace{1mm}

\noindent
{\bf(i)} Si $X'$ est le sous espace de $Y'\times_{Y} V$ constitué des composantes irréductibles sur lesquelles la transformée de $\pi$ (par $\sigma$) est encore génériquement ouverte sur $Y'$, alors, la transformée stricte $\pi':X'\rightarrow Y'$ est ouverte et surjective.\vspace{1mm}

\noindent
{\bf(ii)} L'application induite $\sigma':X'\rightarrow V$ est semi propre et surjective.
\vspace{1mm}

\noindent
{\bf(iii)} $\sigma$ est génériquement ouverte et génériquement finie.\rm
\rm\vspace{2mm}

\noindent
On a  le diagramme
$$\xymatrix{X'\ar@/^1pc/[rr]^{\sigma'}\ar[rd]_{\pi'}\ar@{^{(}->}[r]&Y'\times_{Y} V\ar[d]\ar[r]&X\ar[d]^{\pi}\\
&Y'\ar[r]_{\sigma}&Y}$$

A noter que $\sigma$ n'est ni propre ni même semi propre mais préserve la densité par image inverse. De plus, si ${\rm E}$ est le lieu de dégénérescence de $\pi$ (la réunion des sous ensembles constitué des points en lesquels, les fibres sont de dimension strictement plus grande que ${\rm dim}X -{\rm dim}Y$ cf \lemmaref{lem7'}) et ${\rm N}(Y)$ le lieu non normal de $Y$, $Y'$ s'identifie au produit fibré $Y'\times_{Y} V$ 
en dehors de $\Sigma:=\sigma^{-1}(\overline{\pi({\rm E}\cap V)}\cup {\rm N}(Y))$,
$$X'\setminus \pi'^{-1}(\Sigma)\simeq (Y'\setminus \Sigma)\times_{Y} V$$ 
avec $\sigma: Y'\setminus \Sigma\rightarrow Y$ (resp. $\pi':X'\setminus \pi'^{-1}(\Sigma)\rightarrow X$) ouverte.\vspace{1mm}

\noindent
On voudrait associer à tout morphisme $\pi:X\rightarrow S$ possédant certaines propriétés de bonnes variances, au moins, au niveau des fibres, un morphisme de faisceaux de ${\mathcal O}_{S}$-modules (cohérents si $\pi$ est propre)
$${\mathcal T}^{k}_{\pi,{\mathcal L}}:{\rm I}\!{\rm R}^{n}\pi_{!}{\mathcal L}^{n+k}_{X}\rightarrow{\mathcal L}^{k}_{S} $$
compatible aux restrictions ouvertes sur $X$ et $S$, aux changements de base admissibles en un certain sens et aux morphismes ${\mathfrak T}^{k}_{\pi,\omega}:{\rm I}\!{\rm R}^{n}\pi_{!}{\omega}^{n+k}_{X}\rightarrow{\omega}^{k}_{S} $  du  {\it théorème 1} de \cite{K3}. Selon, les propriétés demandées à $\pi$, on peut proposer différentes constructions.
\subsection{Image directe supérieure pour un morphisme géométriquement plat.} Dans ce cas, on le déduit directement du \theoremref{T0} (dont la preuve se trouve dans \cite{K0}) puisque  $X$ peut-être vu comme le graphe d'une famille  analytique de $n$- cycles $(X_{s})_{s\in S}$ paramétrée par $S$ et pour toute  famille  de supports paracompactifiante   $\Phi$  rencontrant toute famille de supports  $\Psi$  contenant les supports des cycles et telles que $\Phi \cap \Psi$ soit contenue dans la famille des compacts de $X$, permettra de construire 
$$\tilde{\sigma}^{q,0}_{{\Phi},X}:{\rm H}^{n}_{\Phi}(X, {\mathcal 
L}^{n+q}_{Z})\longrightarrow 
{\rm H}^{0}(S, {\mathcal L}^{q}_{S})$$
et donnera naturellement le 
$${\mathcal T}^{k}_{\pi,{\mathcal L}}:{\rm I}\!{\rm R}^{n}\pi_{!}{\mathcal L}^{n+k}_{X}\rightarrow{\mathcal L}^{k}_{S} $$
vérifiant toutes les propriétés fonctorielles annoncées dans le \theoremref{T0}.\vspace{1mm}

\noindent 
Pour donner un caractère intrinsèque au texte, nous rappelons brièvement la construction. Soient $X$ et $S$ deux espaces complexes réduits et $\pi:X\rightarrow S$ un morphisme géométriquement plat à fibres de dimension pure $n$. Pour construire un morphisme de ${\mathcal O}_{S}$-module (cohérents si le morphisme est propre)
$${\mathcal T}^{j}_{\pi}: {\rm I}\!{\rm R}^{n}\pi_{!}{\mathcal L}^{n+j}_{X}\rightarrow{\mathcal L}^{j}_{S}$$ 
( en remplaçant $\pi_{!}$ par $\pi_{*}$ dans le cas propre) compatible aux restrictions ouvertes sur $X$ et aux  changements de base, on procède comme  dans \cite{K0} ou \cite{K1} en commençant par exhiber une telle
flèche relativement à une installation locale de $\pi$ puis
procéder à la globalisation grâce à l'exactitude à droite du foncteur ${\rm I}\!{\rm R}^{n}\pi_{!}$ et à la méthode du découpage.\vspace{2mm}

\indent
 {\bf(a) Construction locale.}\vspace{1mm}

\noindent
Soient  $s_{0}$ un point de $S$ et  $x$ un point de
la fibre $\pi^{-1}(s_{0})$.  Considérons une factorisation locale
de $\pi$ en $x$  donnée (abusivement) par
$\xymatrix{X\ar@/_/[rr]_{\pi}\ar[r]^{f}&Y\ar[r]^{q}&S}$ dans
laquelle $f$ est fini, ouvert  et surjectif sur $Y=S\times U$ avec
$U$ polydisque ouvert relativement compact de ${\Bbb C}^{n}$, $S$ un
ouvert de Stein; $q$ et $p$ étant les projections canoniques sur $S$ et $U$ respectivement, on a, pour tout entier $k$, les décompositions  
$${\mathcal L}^{k}_{Y}=\bigoplus_{i+j=k} q^{*}({\mathcal L}^{i}_{S})\otimes_{{\mathcal O}_{Y}}p^{*}(\Omega^{j}_{U})=\bigoplus_{i+j=k}{\mathcal L}^{i}_{S}\widehat\otimes_{\Bbb C}\Omega^{j}_{U}$$ 
et, donc, une projection naturelle
$${\mathcal L}^{n+q}_{Y}\rightarrow q^{*}({\mathcal L}^{q}_{S})\otimes_{{\mathcal O}_{Y}}
\Omega^{n}_{Y/S}$$
Mais le morphisme canonique de ${\mathcal O}_{S}$-modules
$${{\mathcal L}^{q}_{S}}\otimes_{{\mathcal O}_{S}}{\rm I}\!{\rm R}^{n}q_{!}\Omega^{n}_{Y/S}\rightarrow 
{\rm I}\!{\rm R}^{n}q_{!}(q^{*}({\mathcal L}^{q}_{S})\otimes_{{\mathcal O}_{Y}}
\Omega^{n}_{Y/S})$$
est un isomorphisme puisque, pour tout faisceau cohérent ${\mathcal G}$ sur $S$,  les foncteurs
$${\mathcal G}\rightarrow {\mathcal G}\otimes_{{\mathcal O}_{S}}{\rm I}\!{\rm R}^{n}q_{!}\Omega^{n}_{Y/S},\,\, {\mathcal G}\rightarrow{\rm I}\!{\rm R}^{n}q_{!}(q^{*}({\mathcal G})\otimes_{{\mathcal O}_{Y}}
\Omega^{n}_{Y/S})$$
coincident sur les faisceaux localement libres; le cas général s'en déduit en prenant une résolution localement libre à deux termes pour ${\mathcal G}$. Cela se vérifie aussi rapidement en utilisant une formule de type Künneth sur le produit complété.\vspace{1mm}

\noindent
Alors, en appliquant  le foncteur exact à droite ${\rm I}\!{\rm R}^{n}q_{!}$ à la flèche 
$$f_{*}{\mathcal L}^{n+q}_{X}\rightarrow q^{*}({\mathcal L}^{q}_{S})\otimes_{{\mathcal O}_{Y}}\Omega^{n}_{Y/S}$$ 
on obtient le  morphisme
$${\rm I}\!{\rm R}^{n}q_{!}f_{*}{\mathcal L}^{n+q}_{X}\rightarrow
{\mathcal L}^{q}_{S}\otimes_{{\mathcal O}_{S}}{\rm I}\!{\rm R}^{n}q_{!}\Omega^{n}_{Y/S}$$
Mais, pour tout ouvert
$S'$ de Stein dans $S$, une formule de type K\"unneth donne
$${\rm I}\!{\rm R}^{n}q_{!}\Omega^{n}_{Y/S}(S')\simeq {\mathcal
O}_{S}(S') {{\widehat{\otimes}}{\atop_{{\Bbb C}}}}{\rm H}^{n}_{c}(U,
\Omega^{n}_{U})$$ qui, eu égard à l'isomorphisme ${\rm
H}^{n}_{c}(U, \Omega^{n}_{U})\simeq {\Bbb C}$ donné par la
dualité de Serre, produit, pour $S$ connexe,  l'isomorphisme bien connu pour un morphisme lisse, ${\rm I}\!{\rm
R}^{n}q_{!}\Omega^{n}_{Y/S}\simeq {\mathcal O}_{S}$. 
La finitude de  $f$ garantissant la dégénérescence de la suite spectrale de Leray, et, donc,  
l'isomorphisme de foncteurs ${\rm I}\!{\rm R}^{n}\pi_{!} \simeq {\rm
I}\!{\rm R}^{n}q_{!}f_{*}$,   permet de mettre en évidence la flèche désirée 
$${\rm I}\!{\rm R}^{n}\pi_{!}{\mathcal L}^{n+q}_{X}\rightarrow
{\mathcal L}^{q}_{S}$$
\indent {\bf  (b) Construction globale.}\vspace{1mm}

\noindent
Commençons par remarquer que la construction d'une telle
flèche, qui est toujours de nature locale sur $S$, l'est aussi  sur $X$ puisque la dimension des fibres est inférieure à $n$ et en vertu du  lemme de Reiffen (\cite{R}).  Cela assure, en particulier, 
l'exactitude à droite du foncteur ${\rm I}\!{\rm R}^{n}\pi_{!}$. Par ailleurs, ${\mathcal L}^{n+q}_{X}$ étant de nature locale, on a,  pour tout ouvert  $U$ de $X$ muni de l'injection naturelle $j:U\rightarrow X$, 
$${\mathcal L}^{n+q}_{X}|_{U}={\mathcal L}^{n+q}_{U}$$
Cela étant dit, supposons $X$ paracompact et complètement
paracompact\footnote{Tout ouvert de $X$ est un espace
paracompact; en particulier, pour toute famille paracompactifiante
de supports, l'étendue $\bigcup_{F\in \Phi}F$ est un ouvert
paracompact.}  et considérons un ouvert de Stein (que l'on notera
encore $S$) de $S$ et un recouvrement ouvert localement fini
$(X_{\alpha})_{\alpha\in A}$ de $X$ qui soit $S$-adapté. On peut choisir les ouverts $X_{\alpha}$ de sorte à ce que chacun d'entre eux soit contenu dans une réunion localement finie d'ouverts  munis
d'installations locales du type (avec des notations abusives) 
$$\xymatrix{X_{\alpha}\ar[rdd]_{\pi_{\alpha}}\ar[rr]^{\sigma_{\alpha}}\ar[rd]^{f_{\alpha}}&&Z_{\alpha}\ar[ld]_{p_{\alpha}}\ar[ldd]^{r_{\alpha}}\\
&Y_{\alpha}\ar[d]^{q_{\alpha}}&\\&S&}$$ avec $\sigma_{\alpha}$  un
plongement local, $\pi_{\alpha}$ la restriction de $\pi$ à
$X_{\alpha}$, $Y_{\alpha}$ et $Z_{\alpha}$ sont lisses sur $S$,
$p_{\alpha}$,$q_{\alpha}$ et $r_{\alpha}$ lisses.\vspace{1mm}

\noindent Comme
${\mathcal L}^{n+q}_{X}|_{X_{\alpha}} = {\mathcal L}^{n+q}_{X_{\alpha}}$ on déduit,  du cas local précédent, une collection de morphismes
$${\rm I}\!{\rm R}^{n}{\pi_{\alpha}}_{!}{\mathcal L}^{n+q}_{{X}_{\alpha}}\rightarrow {\mathcal L}^{q}_{S}$$
Mais l'exactitude à droite du foncteur ${\rm I}\!{\rm
R}^{n}\pi_{!}$  qui assure la
surjectivité de la flèche naturelle
$$ \bigoplus_{\atop\alpha\in A}{\rm I}\!{\rm
R}^{n}{\pi_{\alpha}}_{!} {\mathcal L}^{n+q}_{{X_{\alpha}}} \rightarrow
{\rm I}\!{\rm R}^{n}{\pi_{!}}{{\mathcal L}^{n+q}_{X}}$$ permet de
produire, par recollement sur $X$, le morphisme désiré.\vspace{1mm}

\noindent
{\bf(ii)} La compatibilité de la construction aux changements de base est facile à voir car à tout diagramme de changement de base
$$\xymatrix{\widetilde{X}\ar[d]_{\tilde\pi}\ar[r]^{\theta}&X\ar[d]^{\pi}\\
\widetilde{S}\ar[r]_{\nu}&S}$$
est naturellement associé  le diagramme commutatif (qui coincident tous sur la partie régulière de $\pi$)
$$\xymatrix{\nu^{*}{\rm I}\!{\rm R}^{n}{\pi_{!}}{\mathcal L}^{n+q}_{X}\ar[d]_{\nu^{*}({\mathcal T}^{j}_{\pi,{\mathcal L} })}\eq[r]&{\rm I}\!{\rm R}^{n}{\tilde\pi_{!}}\theta^{*}{\mathcal L}^{n+q}_{X}\ar[r]&{\rm I}\!{\rm R}^{n}{\tilde\pi_{!}}{\mathcal L}^{n+q}_{\widetilde X}\ar[d]^{{\mathcal T}^{j}_{\tilde\pi,{\mathcal L}}}\\
\nu^{*}{\mathcal L}^{q}_{S}\ar[rr]&&{\mathcal L}^{q}_{\widetilde S}}$$
\vspace{1mm}

\noindent 
L'isomorphisme de changement de base de la première ligne résulte de :\vspace{1mm}

\indent $\bullet$ l'isomorphisme de changement de base de faisceaux abéliens
$$\nu^{-1}({\rm I}\!{\rm R}^{n}{\pi_{!}}{\mathcal L}^{n+q}_{X})\simeq {\rm I}\!{\rm R}^{n}{\tilde\pi_{!}}\theta^{-1}({\mathcal L}^{n+q}_{X})$$
et de \vspace{1mm}

\indent 
$\bullet$ l'exactitude à droite du foncteur ${\rm I}\!{\rm R}^{n}{\tilde\pi_{!}}$  \vspace{1mm}

\indent qui nous donnent l'isomorphisme
$${\mathcal O}_{\widetilde{S}}\otimes_{\nu^{-1}{\mathcal O}_{S}}\nu^{-1}({\rm I}\!{\rm R}^{n}{\pi_{!}}{\mathcal L}^{n+q}_{X})\simeq {\rm I}\!{\rm R}^{n}{\tilde\pi_{!}}({\mathcal O}_{\widetilde X}\otimes_{\tilde\pi^{-1}{\mathcal O}_{X}}\theta^{-1}{\mathcal L}^{n+q}_{X})$$
c'est-à-dire
$$\nu^{*}({\rm I}\!{\rm R}^{n}{\pi_{!}}{\mathcal L}^{n+q}_{X})\simeq {\rm I}\!{\rm R}^{n}{\tilde\pi_{!}}\theta^{*}({\mathcal L}^{n+q}_{X})$$
La compatibilité aux localisations sur $X$ et $S$ est
aisée puisque c'en est  un cas particulier.\vspace{1mm}

\noindent 
\subsection{Image directe supérieure pour un morphisme propre.}
Grâce au théorème d'aplatissement de Hironaka (\cite{H1}), on peut trouver une modification $\nu:\widetilde{S}\rightarrow S$ et un  diagramme commutatif 
$$\xymatrix{\tilde{X}\ar[rd]_{\psi}\ar[d]_{\tilde\pi}\ar[r]^{\tilde{\nu}}&X\ar[d]^{\pi}\\ 
  \widetilde{S}\ar[r]_{\nu}& S}$$
  dans lequel $\nu$ (resp. $\tilde\nu$) est une modification propre,  $\tilde X$ est la préimage stricte ( qui coincide avec la préimage totale si $\pi$ est ouvert) de $X$ et $\tilde \pi$ un morphisme plat.\vspace{1mm}

  \noindent  Considérons les suites spectrales de second terme
$${\rm E}^{i,j}_{2}:={\rm I}\!{\rm R}^{i}\pi_{!}{\rm I}\!{\rm R}^{j}\tilde{\nu}_{*} {\mathcal L}^{n+k}_{\tilde{X}},\,\,'{\rm E}^{i,j}_{2}:={\rm I}\!{\rm R}^{i}\nu_{*}{\rm I}\!{\rm R}^{j}\tilde{\pi}_{!}{\mathcal L}^{n+k}_{\tilde{X}}$$
 et d'aboutissement  ${\rm I}\!{\rm R}^{i+j}\psi_{*}{\mathcal L}^{n+k}_{\tilde{X}}$ dont on sait qu'elles vérifient  $${\rm E}^{i,j}_{2}=0,\,\forall\,i>n;\,\,'{\rm E}^{i,j}_{2}=0,\,\forall\,j>n$$
 Alors, on a un   morphisme canonique (composition de morphismes latéraux) 
 $${\rm I}\!{\rm R}^{n}{\pi}_{!}{\mathcal L}^{n+k}_{{X}}={\rm I}\!{\rm R}^{n}{\pi}_{!}\tilde{\nu}_{*}{{\mathcal L}}^{n+k}_{\tilde{X}}\rightarrow \nu_{*}{\rm I}\!{\rm R}^{n}\tilde{\pi}_{!}{\mathcal L}^{n+k}_{\tilde{X}}$$
et, par suite, la flèche désirée
$$\xymatrix{{\rm I}\!{\rm R}^{n}{\pi}_{!}{\mathcal L}^{n+k}_{{X}}\ar[r]\ar[d]_{{\mathcal T}^{k}_{\pi,{\mathcal L}}}&\nu_{*}{\rm I}\!{\rm R}^{n}\tilde{\pi}_{!}{\mathcal L}^{n+k}_{\tilde{X}}\ar[d]_{{\mathcal T}^{k}_{\tilde\pi,{\mathcal L}}}\\
 {\mathcal L}^{k}_{S}\ar@{=}[r]&\nu_{*}{\mathcal L}^{k}_{\tilde{S}}}$$
 le diagramme étant  naturellement  commutatif et la compatibilité aux changements de base  découle  de ce qui a été dit  précédemment.\vspace{1mm}

\noindent 
\subsection{Image directe supérieure pour un morphisme génériquement ouvert à fibre générique de dimension $n$.} On utlise, dans ce cas, le théorème d'aplatissement local de \cite{H2} ou de \cite{Si} donnant le diagramme
$$\xymatrix{\widetilde{X}\ar[d]_{\tilde\pi}\ar[r]^{\tilde{\nu}}&X\ar[d]^{\pi}\\ 
  \widetilde{S}\ar[r]_{\nu}& S}$$
  dans lequel $\nu$ (resp. $\tilde\nu$) est un morphisme semi propre,   $\tilde\nu$ un morphisme génériquement ouvert et génériquement fini,  $\tilde X$ est la préimage stricte localisée  de $X$ et $\tilde \pi$ un morphisme équidimensionnel sur $\tilde S$ normal. \vspace{1mm}

  \noindent Pour construire  ${\mathcal T}^{k}_{\tilde\pi,{\mathcal L}}$, il suffit d'utiliser la formule de changement de base précédente associée au produit fibré $\overline X$ et la restriction ${\mathcal L}^{n+k}_{\overline{X}}\rightarrow{\mathcal L}^{n+k}_{\widetilde{X}}$. 

\section{\color{blue}{Image directe supérieure déduite du morphisme ${\mathfrak T}^{k}_{\pi,\omega}:{\rm I}\!{\rm R}^{n}\pi_{!}{\omega}^{n+k}_{X}\rightarrow{\omega}^{k}_{S} $ .}}

Cette approche, quant à elle, de nature plus conceptuelle,  repose sur une construction dont les piliers sont ancrés dans la théorie de la dualité et ne nécessite aucun recours aux théorèmes d'aplatissement comme on l'a vu dans  \cite{K3}.\vspace{1mm}

 \noindent 
En composant avec  l'injection canonique ${\mathcal L}^{n+k}_{X}\rightarrow {\omega}^{n+k}_{X}$, on a, naturellement, le diagramme commutatif
$$\xymatrix{{\rm I}\!{\rm R}^{n}\pi_{!}{\mathcal L}^{n+k}_{X}\ar[rr]\ar[rd]&&{\rm I}\!{\rm R}^{n}\pi_{!}{\omega}^{n+k}_{X}\ar[ld]\\
&{\omega}^{k}_{S}&}$$
Pour lequel, on doit montrer que l'image du morphisme ${\rm I}\!{\rm R}^{n}\pi_{!}{\mathcal L}^{n+k}_{X}\rightarrow{\omega}^{k}_{S}$  est en fait à valeurs dans le faisceau ${\mathcal L}^{k}_{S}$. Pour ce faire, on considère une désingularisation $\nu:\widetilde{S}\rightarrow S$ et le  diagramme de changement de base induit 
$$\xymatrix{\tilde{X}\ar[rd]^{\psi}\ar[d]_{\tilde\pi}\ar[r]^{\tilde{\nu}}&X\ar[d]^{\pi}\\ 
  \widetilde{S}\ar[r]_{\nu}& S}$$
  dans lequel $\phi$ (resp. $\tilde\phi$) est une modification propre et $\pi$ (resp. $\tilde\pi$) est $n$-géométriquement plat. Rappelons que la constante de la dimension des fibres, ouverture et platitude géométrique sont des notions stables par changements de base arbitraires.\vspace{1mm}

  \noindent
  En reprenant les suites spectrales de second terme
$$E^{i,j}_{2}:={\rm I}\!{\rm R}^{i}\pi_{!}{\rm I}\!{\rm R}^{j}\tilde{\nu}_{*} {\mathcal L}^{n+k}_{\tilde{X}},\,\,'E^{i,j}_{2}:={\rm I}\!{\rm R}^{i}\nu_{*}{\rm I}\!{\rm R}^{j}\tilde{\pi}_{!}{\mathcal L}^{n+k}_{\tilde{X}}$$
 et d'aboutissement  ${\rm I}\!{\rm R}^{i+j}\psi_{*}{\mathcal L}^{n+k}_{\tilde{X}}$, on met en évidence, pour les mêmes raisons,  la flèche désirée dans le diagramme commutatif
$$\xymatrix{{\rm I}\!{\rm R}^{n}{\pi}_{!}{\mathcal L}^{n+k}_{{X}}\ar[r]\ar[d]_{{\mathcal T}^{k}_{\pi,{\mathcal L}}}&\nu_{*}{\rm I}\!{\rm R}^{n}\tilde{\pi}_{!}{\mathcal L}^{n+k}_{\tilde{X}}\ar[d]\\
 {\mathcal L}^{k}_{S}\ar@{=}[r]&\nu_{*}\Omega^{k}_{\tilde{S}}}$$
Les propriétés fonctorielles se vérifient sans difficultés comme 
 dans le cas géométriquement plat$\,\blacksquare$

\end{document}